\documentclass[12pt,a4paper]{article}
\newcommand{\blind}{0}
\usepackage{easyReview}
\usepackage{amsmath,amsthm}
\usepackage{amsfonts}
\usepackage{amssymb}
\usepackage{bbm}
\usepackage{mathtools}
\usepackage{pdflscape}

\usepackage{url}

\usepackage[version=4]{mhchem}
\usepackage{stmaryrd}

\usepackage{subcaption} 
\usepackage{graphicx}
\usepackage{pgfplots}
\usepackage[all]{nowidow}
\usepackage[utf8]{inputenc}
\usepackage{tikz}
\usepackage{multicol}
\usepackage{algpseudocode,algorithm,algorithmicx}
\usepackage{scalerel}
\usepackage{natbib}
\usepackage{bibunits}
\usepackage[normalem]{ulem}

\usepackage[english]{babel}

\usepackage{mathdots}
\usepackage{yhmath}
\usepackage{cancel}
\usepackage{color}
\usepackage{array}
\usepackage{multirow}
\usepackage{gensymb}
\usepackage{tabularx}
\usepackage{booktabs}
\usetikzlibrary{fadings}
\usetikzlibrary{patterns}
\usetikzlibrary{shadows.blur}

\usepackage{enumerate}
\usepackage[algo2e]{algorithm2e}

\newtheorem{theorem}{Theorem}
\newtheorem{lemma}{Lemma}

\newtheorem{assumption}{Assumption}

\newtheorem{remark}{Remark}
\newtheorem{definition}{Definition}

\newcommand{\Z}{\mathbb{Z}}

\newcommand{\R}{\mathbb{R}}

\newcommand{\cK}{\mathcal{K}}

\newcommand{\bH}{\mathbb{H}}

\newcommand{\wt}{\widetilde}

\newcommand{\X}{\mathcal{X}}

\newcommand{\N}{\mathcal{N}}

\newcommand{\cH}{\mathcal{H}}

\newcommand{\F}{\mathcal{F}}

\newcommand{\li}{\langle}
\newcommand{\ri}{\rangle}
\newcommand{\fr}{{\lfloor rn\rfloor}}
\newcommand{\ftr}{{\lfloor \wt rn\rfloor}}
\newcommand{\frr}{{\lfloor r_0n\rfloor}}

\newcommand{\fb}{{\lfloor bn \rfloor}}

\DeclareMathOperator*{\var}{Var}
\DeclareMathOperator*{\cov}{Cov}

\DeclareMathOperator{\E}{\mathbb{E}}
\newcommand{\id}{\mathbbm{1}}

\usepackage[nodisplayskipstretch]{setspace}
\begin{document}
	\begin{bibunit}[chicago]
	
	\def\spacingset#1{\renewcommand{\baselinestretch}%
		{#1}\small\normalsize} \spacingset{1}

	
	\if0\blind
	{
		\title{\bf Online Change-Point Monitoring for Object-valued Time Series}
		\author{Yi Zhang, Tim Kutta and Xiaofeng Shao
			\thanks{Yi Zhang (corresponding author) is a postdoctoral fellow, and Xiaofeng Shao is a professor in the Department of Statistics and Data Science at Washington University in St Louis. Emails: {\tt yiz@wustl.edu} and {\tt shaox@wustl.edu}. Tim Kutta is an assistant professor in the Department of Mathematics at Aarhus University, Denmark, and is partially supported by NNF grant 52553. Shao's research is partially supported by NSF grant DMS-2524356.}\\
		}
		\maketitle
	} \fi
	
	\if1\blind
	{
		\bigskip
		\bigskip
		\bigskip
		\begin{center}
			{\LARGE\bf Title}
		\end{center}
		\medskip
	} \fi
	
	\bigskip
	\begin{abstract}
		We develop closed- and open-end procedures for monitoring changes in the marginal distribution of object-valued time series. The method combines two distance-based Hilbert-space embeddings, a monitoring-time-dependent projection, and self-normalization. It is computable entirely from pairwise distances, does not require long-run variance estimation, and admits exact recursive updates. Under weak temporal dependence, the closed-end null limit is a pivotal Brownian functional. For monitoring over an unbounded horizon, we introduce a growing monitoring-time weight, establish a maximal inequality and uniform remote-tail control, and derive an open-end pivotal limit together with an equivalent fixed-interval representation for critical-value simulation. Both procedures are consistent against fixed marginal changes and retain the $n^{-1/2}$ and $n^{-1/4}$ local detection boundaries associated with the linear and quadratic projection signals. Simulations with distribution-valued time series illustrate favorable numerical properties, and an application to monthly stock-return distributions illustrates the usefulness of our procedure. Further simulations with graph-valued time series and an application to spatial point processes of seismic activity are reported in the supplementary material.
		
	\end{abstract}
	
	\noindent%
	\vfill
	
	\newpage
	
	\spacingset{1.667} 
	
	\section{Introduction}
	
	Modern data streams increasingly consist of complex objects rather than scalars or Euclidean vectors. Examples include sequences of probability distributions, covariance matrices, shapes, images, and networks. A useful abstraction treats each observation as a random element of a metric space, thereby retaining the geometry relevant to the application without imposing an artificial coordinate representation. This generality comes with a fundamental difficulty: addition, scalar multiplication, inner products, and other operations underlying classical CUSUM statistics may not be defined, and the information available to the analyst may consist only of pairwise distances between observations. 	These issues are particularly consequential in online monitoring: after an initial training sample serving as a benchmark, observations arrive sequentially, and the statistical task is to detect a change in the data-generating mechanism in real time. For an introduction to the monitoring problem, we refer to the review by \citet{aue:kirch:2024}. Developing a monitoring procedure that is broadly sensitive to changes in the marginal distribution, accommodates serial dependence, and remains implementable for general random objects is an interesting but nontrivial problem.	
	
	The literature on change-point inference for object-valued data has developed mainly in the retrospective (offline) setting. Graph-based scan statistics accommodate multivariate and non-Euclidean observations through similarity graphs, but their theory  primarily covers independent data \citep{chen2015graph,chu2019asymptotic}. Fr\'echet means and variances provide notions of location and dispersion for random objects in metric spaces \citep{frechet1948,dubey2019frechet}. Change-point procedures based on these quantities test for changes in the corresponding distributional features \citep{dubey2020frechet,jiang2024two}, whereas distance-profile methods can target broader distributional changes under independence \citep{dubey2026distanceprofiles}. Most closely related to the present work, \citet{zhang2025change} combine Hilbert-space embedding, sample splitting, projection, and self-normalization to construct a retrospective omnibus test for the marginal distribution of a weakly dependent object-valued time series. Their statistic has a pivotal null limit and can be evaluated from pairwise distances, but the sample-splitting device used to estimate the change direction is inherently retrospective.
	
	Sequential monitoring has a separate and extensive history, beginning with procedures based on a stable historical period and continuing through methods for general finite-dimensional parameters and high-dimensional data \citep{chu1996monitoring,Holger2020,wu2022adaptive}. Omnibus distributional monitoring is available for multivariate Euclidean time series through empirical distribution functions, with data-dependent thresholds obtained by resampling or covariance estimation \citep{kojadinovic2021nonparametric,holmes2024multipurpose}. Recent work has also considered monitoring mean changes in functional or Banach-space-valued observations \citep{kutta2025monitoring,zhu2026online} and changes in specialized streams of distributions \citep{horvath2021monitoring,zeng2026beyond}. A contemporaneous approach based on degenerate $U$-statistics treats distributional changes for observations in a separable metric space under independence; its null calibration depends on the spectrum of a kernel operator \citep{boniece2026sequential}. These developments leave an important intersection unresolved:  omnibus consistent monitoring for changes in the entire marginal distribution of general, dependent object-valued time series. Our aim is more ambitious still. We seek a procedure that is essentially tuning-free, requiring neither data-dependent resampling nor bandwidth selection for long-run covariance estimation. The procedure should also depend only on pairwise distances between the observed objects, thereby avoiding another potentially consequential modeling choice: the construction of an explicit feature representation-through coordinates, basis expansions, embeddings, or other problem-specific transformations-on which the subsequent analysis would depend. A sample-splitting strategy along the lines of \citet{zhang2025change} faces a fundamental obstacle in sequential monitoring: its sequential extension does not  deliver a pivotal limiting distribution. Section~\ref{sec2_2} explains why, and the supplementary material provides the corresponding mathematical details.

	For these reasons, we propose a new monitoring approach, where a time-dependent projection direction is formed at each step from all currently available data. Because this projection is based on a mean rather than a difference of means, using only a single (first) Hilbert-space embedding would make its direction dependent on the arbitrary origin of that embedding and could make its population limit vanish. We therefore introduce a second embedding through an integrally strictly positive definite radial kernel. The resulting feature representation is bounded, has a nonzero population mean, and has all required inner products determined by pairwise distances; it also distinguishes any two distinct marginal distributions. The projected process and the historical self-normalizer (i.e., the self-normalizer computed based on the historical stable sample) then share the same long-run dependence structure, which cancels to yield a pivotal Brownian-functional null limit under weak temporal dependence.
	
	The proposed construction has six main features. First, it is asymptotically omnibus for fixed changes in the marginal distribution within the stated negative-type metric-space framework. The second embedding is based on an integrally strictly positive definite kernel and therefore makes distinct marginal distributions have distinct RKHS means. 
	The omnibus property is non-trivial, because the projection direction of our statistic can be orthogonal to the change direction at a certain point in time. However, by construction, the projection direction evolves over time, and in fact rotates in the direction of the change. This leads to asymptotic power against arbitrary fixed alternatives.
	Theorem~\ref{the_power} further quantifies local sensitivity: the detection boundary is $n^{-1/2}$ when the projection has a nonzero first-order alignment with the mean change. The rotation of the projection direction then produces an additional quadratic signal and an $n^{-1/4}$ boundary when that alignment vanishes.  Section~\ref{sec_power_enhancement} discusses a power-enhancement method, inspired by \citet{fan2015power}, that recovers detection rates arbitrarily close to $n^{-1/2}$ even in the orthogonal case.
	Second, under weak serial dependence, Theorem~\ref{the_null} establishes a pivotal closed-end limiting null distribution expressed solely through a standard Brownian motion. Third, Section~\ref{sec_open} extends the method to an unbounded monitoring horizon; Theorems~\ref{the_open_null}--\ref{the_open_power} establish a pivotal open-end null limit and retain both local detection boundaries for fixed scaled change locations. Fourth, the historical self-normalizer removes the need to estimate a long-run variance or covariance operator and therefore avoids bandwidth, block-length, and truncation choices associated with dependence estimation. Fifth, all inner-product evaluations required by the statistic reduce to transformed pairwise distances, so neither Hilbert-space embedding needs to be constructed explicitly. Finally, the statistic admits an exact recursive implementation. Each new observation requires $O(k)$ distance evaluations and arithmetic operations at monitoring time $k$, giving $O(n^2)$ total computation through a fixed monitoring horizon and $O(n)$ working storage apart from the observations, compared with $O(n^3)$ operations for direct computation.
	
	The rest of the paper is organized as follows. Section~\ref{sec2single} develops the closed-end statistic, its null and alternative theory, its recursive implementation and a power enhancement method. Section~\ref{sec_open} introduces the open-end weight function, establishes the open-end null and local-power theory, and gives a fixed-interval method for simulating critical values. Section~\ref{sec3simu} reports a simulation study for distribution-valued time series. Section~\ref{sec4real} presents an application to monthly stock-return distributions. Section~\ref{sec5conclu} summarizes our main conclusions. The supplementary material contains a simulation study for graph-valued time series, a second data application to micro-earthquakes, and proofs of the mathematical results. The supplementary material is not included in this arXiv version.
	
	The following notation will be used throughout the paper. For a positive integer $n$, define $[n]=\{1,2,\dots,n\}$. Let $\|\cdot\|$ be the Euclidean norm on $\R^p$ and $\|\cdot\|_\cH=\li\cdot,\cdot\ri_{\cH}^{1/2}$ be the norm on the Hilbert space $\cH$. For $\R$- or $\cH$-valued random variables $\{X_t\}_{t\in\Z}$, let
	$\F_a^b=\sigma(X_t:a\leq t\leq b)$ for any
	$-\infty\leq a\leq b\leq\infty$. For a fixed $T>1$, define the triangular index set $\Delta_T=\{(\tau,r):1\leq r\leq\tau\leq T\}$. For a compact interval $I\subset\R$ and a normed space $E$, let $D_E(I)$ denote the space of $E$-valued cadlag functions on $I$ endowed with the Skorokhod $J_1$ metric \citep[Section~12]{billingsley1999convergence};
	we write $D_E[a,b]$ when $I=[a,b]$. For an index set $A$, let $\ell^\infty(A,E)$ denote the space of bounded functions $g:A\to E$ equipped with the supremum norm $\|g\|_\infty=\sup_{a\in A}\|g(a)\|_E$, and write $\ell^\infty(A)=\ell^\infty(A,\R)$ \citep[Section~1.5]{vandervaart1996weak}. Thus, $\ell^\infty(\Delta_T)$ is the space of bounded real-valued functions on $\Delta_T$. We use ``$\rightsquigarrow$'' to denote weak convergence in the function space specified in context.

	\section{Closed-End Monitoring: Methodology and Theory}\label{sec2single}
	
	Before introducing the monitoring problem, we first collect the probabilistic preliminaries needed for the asymptotic analysis. Section~\ref{sec2_1} states the relevant definitions and a functional central limit theorem for weakly dependent Hilbert-space-valued time series. Section~\ref{sec2_2} then formulates the closed-end monitoring problem and constructs the proposed statistic. Section~\ref{sec_closed_theory} develops the asymptotic results under null and alternative, while Section~\ref{sec_iterative} presents the recursive implementation method.
	
	\subsection{Preliminary Result}\label{sec2_1}

	Given an $\R$- or $\cH$-valued process
	$\{X_t\}_{t\in\Z}$, write
	$\F_I=\sigma(X_t:t\in I)$ for $I\subset\Z$. For a positive
	integer $m$, define the interlaced $\rho$-mixing coefficient used
	below \citep{bradley2005,BUCC2017} by
	\[
	\rho(m)=\sup_{\substack{I,J\subset\Z,\ I,J\neq\varnothing\\
			\operatorname{dist}(I,J)\geq m}}
	\rho(\F_I,\F_J),
	\]
	where
	$\operatorname{dist}(I,J)=\min\{|i-j|:i\in I,j\in J\}$ and
	\[
	\rho(\mathcal A,\mathcal B)
	=\sup\left\{
	\frac{|\cov(U,V)|}
	{\sqrt{\var(U)\var(V)}}:
	U\in L_2(\mathcal A),\ V\in L_2(\mathcal B),\
	\var(U),\var(V)>0
	\right\}.
	\]
	Here, $L_2(\mathcal A)$ is the set of real-valued,
	$\mathcal A$-measurable random variables with finite second
	moment. The process is called $\rho$-mixing if
	$\rho(m)\to0$ as $m\to\infty$. As in \cite{BUCC2017}, define
	$\alpha_{1,1}(m)=\sup_{j\in\Z,k\geq m}
	\alpha(\F_j^j,\F_{j+k}^{j+k})$, where
	\[
	\alpha(\mathcal A,\mathcal B)
	=\sup\{|P(A\cap B)-P(A)P(B)|:
	A\in\mathcal A,\ B\in\mathcal B\}.
	\]
	A linear operator $Q:\cH\to\cH$ is said to be a positive, self-adjoint, trace class operator if (i) $\li Qh,h\ri_\cH\geq0$ for any $h\in\cH$;
	(ii) $\li Qh,x\ri_\cH=\li h, Qx\ri_\cH$ for any $h,x\in\cH$; (iii) $\sum_{k=1}^{\infty}\li Qe_k,e_k \ri_\cH <\infty$ for any orthonormal base $\{e_k\}$ of $\cH$.
	
	Based on these concepts, we can define $\cH$-valued Gaussian random element as well as $\cH$-valued Brownian motion.
	\begin{definition}[\cite{andresen2013}]
		An $\cH$-valued random element $X$ is Gaussian with covariance operator $Q$ and mean $\mu\in \cH$ (denoted as $X\sim \N(\mu,Q)$) if for any $h\in \cH$, $\li X,h\ri_\cH \sim \N(\li \mu,h\ri_\cH, \li Qh,h\ri_\cH)$. 
	\end{definition}
	
	\begin{definition}[\cite{andresen2013}] \label{def2}
		Let $T>0$ be fixed and let $Q$ be a positive, self-adjoint, trace-class operator on $\cH$. An $\cH$-valued Brownian motion with covariance operator $Q$ on $[0,T]$ is a stochastic process $\{B_Q(t)\}_{t\in[0,T]}$ such that (i) $B_Q(0)=0$; (ii) $B_Q(t)$ has continuous sample paths; (iii) $B_Q(t)$ has independent increments; and (iv) $B_Q(t)-B_Q(s)\sim \N(0,(t-s)Q)$ for any $0\leq s<t\leq T$.
		
	\end{definition}

	The following FCLT for stationary $\cH$-valued random variables is from \cite{BUCC2017}.
	
	\begin{theorem}\label{th1}
		Let $\{X_t\}_{t\in \Z}$ be a strictly stationary $\cH$-valued random process with mean $\mu$. Assume $\{X_t\}_{t\in \Z}$ is $\rho$-mixing and the following conditions hold for some $\delta>0$:
		\begin{enumerate}
			\item \label{th1_1} $\E\|X_1\|_\cH^{2+\delta}<\infty$,
			\item \label{th1_2} $\sum_{n=1}^{\infty}[\alpha_{1,1}(n)]^{\delta/(2+\delta)}<\infty$.
		\end{enumerate}
		Then, for every fixed $T>0$,
		\begin{align*}
			\left\{\frac{1}{\sqrt{n}}\sum_{t=1}^{\lfloor nr\rfloor}\big(X_t-\mu\big)\right\}_{r\in[0,T]}
			\rightsquigarrow
			\{B_Q(r)\}_{r\in[0,T]}
			\quad\text{in }D_\cH[0,T],
		\end{align*}
		where $\{B_Q(r)\}_{r\in[0,T]}$ is a Brownian motion with
		\begin{align}\label{eq11}
			\li Qx, y\ri_\cH = \sum_{t\in\Z}\E\big(\li X_0-\mu,x \ri_\cH\li X_t-\mu,y \ri_\cH\big) \mbox{ for any } x,y\in\cH.
		\end{align}
		Furthermore, the series in Equation (\ref{eq11}) converges absolutely.
	\end{theorem}

	\subsection{Problem Setup and Closed-End Monitoring Statistic}\label{sec2_2}
	
	From now on, we assume the time series $\{X_t\}_{t\in\Z}$ takes values in a separable metric space $(\X,d)$ of negative type (see Definition 2.2 in \cite{chak2021} or \cite{lyons2013}). Let $P_i$ be the marginal distribution of $X_i$. In this section, we are concerned with the  closed-end formulation of the monitoring problem, as discussed e.g. in \cite{Holger2020} and \cite{wu2022adaptive}.
	We assume that initially we observe a historical sample of $n$ random objects $\{X_t\}_{t=1}^n$ 
	and afterwards start the  monitoring procedure, where data arrive sequentially, first $X_{n+1}$, then $X_{n+2}$ and so on. The aim is to detect a change in the marginal distribution, and at each
	time point $k>n$ we decide whether a change has occurred. In the
	closed-end setting, we consider monitoring times
	$k=n{+}1,\dots,Tn$, where $T>1$ is a constant. Our approach is viable for any $T>1$, but in our subsequent presentation, we mostly assume that $T$ is an 
	integer for simplicity. Monitoring ends at time $Tn$ regardless of
	whether a change point is detected. Under the null hypothesis, we have $P_i=P_j$ for any $i,j\in[Tn]$. Under the alternative, for some fixed $r_0\in(1,T)$ and $k^\ast=\frr$, we assume $P_1=\cdots=P_{k^\ast}\not=P_{k^\ast+1}=\cdots=P_{Tn}$.
	
	For the offline change point detection problem, where we only observe $\{X_t\}_{t=1}^n$, \cite{zhang2025change} proposed a test procedure which combines Hilbert space embedding, sample splitting and self normalization (SN). Specifically, by Proposition 3.2 in \cite{chak2021}, 
	for a metric space $(\X,d)$ of negative type, there exists a Hilbert space $\cH$ and an embedding map $\phi:\X\rightarrow \cH$ such that $d(x,x')=\|\phi(x)-\phi(x')\|_{\cH}^2$. For $t\in \Z$, let $Y_t=\phi(X_t)$. If $\int\|Y_t\|_\cH dP_t<\infty$, define $\xi_t = \E[ Y_t]$, which is the unique element in $\cH$ corresponding to the bounded linear functional $\cH\to\R: h\to \E \li h,Y_t \ri_\cH $. If we further assume the metric space $(\X,d)$ to be of strong negative type \citep[][Definition 2.2]{chak2021}, which is assumed in \cite{zhang2025change},  by Proposition 3.1 in \cite{lyons2013}, the null hypothesis is equivalent to $H_0:\xi_1=\cdots=\xi_{Tn}$
	versus
	$H_1:\xi_1=\cdots=\xi_{k^\ast}\not=\xi_{k^\ast+1}=\cdots=\xi_{Tn}$.
	
	Since the SN approach cannot be directly applied for infinite dimensional (Hilbert space valued) time series \citep[][Section 2.1]{zhang2025hypothesis}, for a fixed $b\in(0,1/2)$ and $\fb=\lfloor bn\rfloor$, \cite{zhang2025change} first split the data into three parts: $\X_1=\{X_1,\dots,X_{\fb}\}$, $\X_2 =\{X_{\fb+1},\dots,X_{n-\fb}\} $ and $\X_3 = \{X_{n-\fb+1},\dots,X_n\}$.
	\begin{samepage}
		Based on $\X_1$ and $\X_3$, they defined the projection direction
		\begin{align}\label{eq_proj_zhang}
			\bar Y_1-\bar Y_3=\frac{1}{\sqrt{n}} \sum_{t=1}^{\fb} Y_t-\frac{1}{\sqrt{n}} \sum_{t=n-\fb+1}^{n} Y_t,
		\end{align}
		which estimates the mean change $\xi_1{-}\xi_n$ (subject to a rescaling factor) under $H_1$.
	\end{samepage}
	They then projected the embedded observations in $\X_2$ along the direction $\bar{Y}_1-\bar{Y}_3$ by defining
	$Z_t=\li\bar{Y}_1-\bar{Y}_3,Y_t\ri_\cH$ for $t=\fb{+}1,\cdots,n{-}\fb$. Finally, they applied the SN mean-change statistic of \citet{shaozhang2010} to $\{Z_t\}$. By \citet[Theorem~2(i)]{zhang2025change}, the resulting statistic has a pivotal limiting null distribution.
	
	However, the idea of sample splitting does not work in the online setting. The main reason is that more and more data are observed as the monitoring time $k$ progresses and the third subsample $\X_3$ in \cite{zhang2025change}, which corresponds to $\{X_{k-\fb+1},\dots,X_k\}$ in the online setting, will change with $k$. As a consequence, the limiting null distribution for the statistic defined following the sample splitting idea in \cite{zhang2025change} is not pivotal; see the supplementary material for more details.
	
	To tackle this problem, we construct a different projection direction using the average of all observations available at each monitoring time after a novel embedding scheme. A direct average after the embedding $\phi$ is not suitable: because $\phi$ is determined by $d$ only up to a translation, such an average is neither an intrinsic distance-based direction nor guaranteed to have a nonzero population limit. We therefore further embed $Y_t$ in a reproducing kernel Hilbert space (RKHS) on $\cH$. Specifically, define the continuous, positive-definite kernel $\cK:\cH\times \cH\to \R$ as $\cK(a,b) = f(\|a-b\|^2_\cH)$ where $f:[0,\infty)\to \R$ is defined as
	\begin{align}
		f(x) = \int_{[0,\infty)}e^{-tx}dv(t),
	\end{align}
	with $v(\cdot)$ being any finite Borel measure satisfying $v((0,\infty))>0$. According to Corollary 3.2 in \cite{ziegel2024characteristic}, $\cK$ is integrally strictly positive definite with respect to all finite signed Borel measures on $\cH$; equivalently, every nonzero such measure has a nonzero RKHS embedding. For a probability measure $R$ on $\X$, write
	\begin{align*}
		\mu_R=\int_{\X}\cK(\phi(x),\cdot)\,dR(x)
	\end{align*}
	for its mean embedding in the RKHS (denoted as $\bH$). Since $\phi(\cdot)$ is injective, Theorem 3.9 in \cite{parthasarathy2005probability} ensures that distinct probability measures on $\X$ induce distinct probability measures on $\cH$. Integral strict positive definiteness therefore gives
	\begin{equation}\label{eq_second_embedding_properties}
		\begin{aligned}
			\text{Nonzero embedding:}\quad&
			\|\mu_R\|_\bH^2
			=\int_{\X}\int_{\X}f\{d(x,x')\}\,dR(x)\,dR(x')>0,\\
			\text{Injectivity:}\quad&
			R\neq S
			\ \Longrightarrow\ \|\mu_R-\mu_S\|_\bH^2>0,\\
			\text{Nonzero mixtures:}\quad&
			(1-a)\mu_R+a\mu_S
			=\mu_{(1-a)R+aS}\neq0,\qquad a\in[0,1].
		\end{aligned}
	\end{equation}
	\paragraph{Role of the second embedding.}
	The nonzero-embedding and nonzero-mixture properties in Equation (\ref{eq_second_embedding_properties}) ensure that, under the null, the projection direction $\bar\phi(k)$ defined in Equation (\ref{eq_noisepro}) converges to the nonzero direction $\mu_{P_1}$ and, under the alternative, its limit at every monitoring time remains the nonzero mean embedding of a mixture of the pre- and post-change distributions. The injectivity property in Equation (\ref{eq_second_embedding_properties}) further gives $P_1\neq P_{Tn}\Rightarrow\mu_{P_1}\neq\mu_{P_{Tn}}$. Thus, every marginal distributional change becomes a nonzero RKHS mean change, which is the basis of the omnibus-consistency argument below. Denote $\wt Y_t=\cK(Y_t,\cdot)$ and $\mu_t=\mu_{P_t}=\E\wt Y_t$. Then the null and alternative hypotheses are equivalent to $H_0:\mu_1=\cdots=\mu_{Tn}$
	versus
	$H_1:\mu_1=\cdots=\mu_{k^\ast}\not=\mu_{k^\ast+1}=\cdots=\mu_{Tn}$.
	
	For $t\geq1$, define the partial-sum process $U_t=\sum_{i=1}^{t}\wt Y_i$. At a given time $k\geq n{+}2$, we use the sample average of the second-embedded features as the projection direction $\bar\phi(k)$, defined by
	\begin{align}\label{eq_noisepro}
		\bar \phi(k) = \frac{U_k}{k}.
	\end{align}
	For $s\in\{n+1,\dots,k-1\}$, define the Hilbert-space CUSUM $C_{k,n}(s)$ by
	\begin{align}\label{eq_cusum_process}
		C_{k,n}(s)=\frac{1}{\sqrt{n}}\left(U_s-\frac{s}{k}U_k\right)
	\end{align}
	and define its scalar projection $T_{k,n}(s)$ by
	\begin{align}\label{eq_tkn}
		T_{k,n}(s)=\li\bar\phi(k),C_{k,n}(s)\ri_\bH
		=\frac{\li U_k,U_s\ri_\bH}{k\sqrt{n}}
		-\frac{s\|U_k\|_\bH^2}{k^2\sqrt{n}}.
	\end{align}
	The historical self-normalizer $W_n$ is constructed once from the historical sample as follows:
	\begin{align}\label{eq_wn}
		W_n
		=&\frac{1}{n^4}\sum_{t=1}^{n}
		\left\langle U_n ,U_t-\frac{t}{n}U_n\right\rangle_\bH^2.
	\end{align}
	Our monitoring statistic at time $k$ is defined as 
	\begin{align} \label{eq_gn_size}
		G_n(k) =	 \max_{s=n+1,n+2,\dots,k-1}\frac{|T_{k,n}(s)|}{\sqrt{W_{n}}}.
	\end{align}
	Define $G_n = 	\max_{k = n+2,\dots,Tn} G_n(k)$. As shown in Theorem \ref{the_null} below, $G_n$ converges to a pivotal limiting distribution $G$ under $H_0$. For some $\gamma\in(0,1)$, let $G_\gamma$ be the $100(1-\gamma)$th upper percentile of $G$, then an asymptotic level-$\gamma$ monitoring scheme is that: at each monitoring time $k$, we declare that  a change has occurred and stop if $G_n(k)\geq G_\gamma$. Note that a similar idea of constructing the self-normalizer solely from the historical sample was used by \citet{chan2021self} for Euclidean-valued time series. Since $W_n$ remains fixed throughout the monitoring period, it needs to be computed only once, which reduces the cost of sequential updating; see Section~\ref{sec_iterative}.
	
	\paragraph{Omnibus consistency.}
	Equation (\ref{eq_second_embedding_properties}) ensures that the population projection direction is nonzero at every monitoring fraction. Nevertheless, at a particular monitoring time, this direction could be orthogonal to the change direction. The maximization over monitoring times prevents such an isolated cancellation from eliminating power. To see this, consider the fixed alternative with $k^\ast=\frr$ and define $\Delta_\mu=\mu_{Tn}-\mu_1\neq0$, where the inequality again follows from Equation (\ref{eq_second_embedding_properties}). For any fixed $r\in(r_0,T)$, let $k=\fr$ and $a=1-r_0/r\in(0,1-r_0/T)$. Under Assumption \ref{assump1}, the laws of large numbers yield
	\begin{align} \label{e:drift}
		&\bar\phi(k)
		\stackrel{p}{\to}(1-a)\mu_1+a\mu_{Tn}\quad \text{ and }\quad \frac{1}{\sqrt n}C_{k,n}(k^\ast)
		\stackrel{p}{\to}-r_0a\Delta_\mu.
	\end{align}
	Consequently,
	\begin{align}
		&\frac{1}{\sqrt n}T_{k,n}(k^\ast)
		\stackrel{p}{\to}-r_0a\,h(a),\label{eq_power_signal} \qquad \textnormal{where}\\
		& h(a)
		=\li(1-a)\mu_1+a\mu_{Tn},\Delta_\mu\ri_\bH=\li\mu_1,\Delta_\mu\ri_\bH+a\|\Delta_\mu\|_\bH^2.\nonumber
	\end{align}
	Since $h(a)$ is affine in $a$ with strictly positive slope $\|\Delta_\mu\|_\bH^2$, it can vanish at no more than one value of $a$. We may therefore choose a fixed $\wt r\in(r_0,T)$, with $\wt a=1-r_0/\wt r$, such that $h(\wt a)\neq0$. For $\wt k=\ftr$, Equation (\ref{eq_power_signal}) then gives
	\begin{align*}
		|T_{\wt k,n}(k^\ast)|
		=\sqrt n\,r_0\wt a\,|h(\wt a)|+o_p(\sqrt n).
	\end{align*}
	The historical self-normalizer is computed entirely from the stable pre-change sample, and Assumption \ref{assump1} also gives $W_n=O_p(1)$. It follows that $	G_n
		\geq\frac{|T_{\wt k,n}(k^\ast)|}{\sqrt{W_n}}
		\stackrel{p}{\longrightarrow}\infty.$
	Thus, for every $P_1\neq P_{Tn}$, the rejection probability at any fixed critical value converges to one. A particular projection direction may miss the change at one monitoring fraction, but scanning over $k$ ensures that the procedure retains asymptotic power against every marginal distributional change.

	\begin{remark}
		Although the feature sequence $\{\wt Y_t\}$ is defined through the embedding map $\phi(\cdot)$, this map may be unknown or difficult to construct explicitly. Nevertheless, for any positive integers $a$ and $b$,
		\begin{align}
			\li U_a,U_b\ri_\bH
			=&\sum_{i=1}^{a}\sum_{j=1}^{b}\li\wt Y_i,\wt Y_j\ri_\bH
			=\sum_{i=1}^{a}\sum_{j=1}^{b}f(d(X_i,X_j)).
		\end{align}
		Therefore, by Equations (\ref{eq_tkn}) and (\ref{eq_wn}), $G_n(k)$ is determined entirely by the pairwise distances among the observations available through time $k$, and neither embedding needs to be constructed explicitly.
		
		This distance-only computability is an additional role of the second embedding. If only the first embedding $Y_t=\phi(X_t)$ were used with the same sample-average projection, then the analogous statistic would require knowledge of $\phi$ for its calculation. Specifically, let $\bar Y(k)=k^{-1}\sum_{i=1}^kY_i$. The statistic corresponding to $T_{k,n}(s)$ would be 
		\begin{align*}
			T_{k,n}^\phi(s)
			={}&\frac{1}{\sqrt n}\left\{
			\sum_{t=1}^s\langle\bar Y(k),Y_t\rangle_\cH
			-\frac{s}{k}\sum_{t=1}^k\langle\bar Y(k),Y_t\rangle_\cH
			\right\}\\
			={}&-\frac{1}{2k\sqrt n}\left\{
			\sum_{i=1}^k\sum_{t=1}^s d(X_i,X_t)
			-\frac{s}{k}\sum_{i=1}^k\sum_{t=1}^k d(X_i,X_t)
			\right\}
			+R_{k,n}^\phi(s),
		\end{align*}
		where $	R_{k,n}^\phi(s)
		=\frac{1}{2\sqrt n}\left\{
		\sum_{t=1}^s\|\phi(X_t)\|_\cH^2
		-\frac{s}{k}\sum_{t=1}^k\|\phi(X_t)\|_\cH^2
		\right\}$ depends on $\phi$. Calculation of this average-projection statistic would therefore require knowledge of $\phi$. The second embedding removes this dependence because every required inner product equals $f(d(X_i,X_j))$.
	\end{remark}
	
	\begin{remark}\label{rmk_2}
		In implementation, a potential concern for statistics based on kernel embeddings is the disparity between diagonal and off-diagonal inner products. For example, consider $f(x)=\exp(-x)$, for which $\li\wt Y_i,\wt Y_j\ri_\bH=\exp\{-d(X_i,X_j)\}$. In this case, $\|\wt Y_i\|_\bH^2=f(0)=1$, whereas $\li\wt Y_i,\wt Y_j\ri_\bH=f\{d(X_i,X_j)\}$ can be close to zero when $d(X_i,X_j)$ is large. Thus, the diagonal terms may appear prominent in sums of pairwise inner products, although they are identical for all observations and contain no information about a distributional change.
		
		For the proposed monitoring statistic, however, this is less of an issue because all terms of the form $\|\wt Y_i\|_\bH^2$ cancel from both the numerator and denominator. Specifically, these terms contribute $sf(0)/(k\sqrt n)$ to both terms in the numerator defined in Equation~(\ref{eq_tkn}) and therefore cancel. The same cancellation occurs in $W_n$. Consequently, the numerator and denominator of the proposed statistic are roughly of the same scale (i.e., both are of the same order as $\li\wt Y_i,\wt Y_j\ri_\bH$ for $i\neq j$). 
	\end{remark}

	\subsection{Theoretical Properties under the Null and Alternative}\label{sec_closed_theory}
	
	We begin with a formal assumption describing weak dependence and stationarity before and, if a change occurs, after the change. Following \cite{kutta2025monitoring}, we describe pre- and post-change distribution together, using the jointly stationary $\bH\times\bH$-valued time series
	$\{(\widetilde Y_t^{(1)},\widetilde Y_t^{(2)})\}_{t\in\mathbb Z}$. The first component represents the pre- and the second component the post-change regime. For a change point $k^\ast$, define the observed process by
	\begin{align}\label{e:Y:case}
		\widetilde Y_t :=
		\begin{cases}
			\widetilde Y_t^{(1)}, & 1\leq t\leq k^\ast,\\
			\widetilde Y_t^{(2)}, & k^\ast<t\leq Tn.
		\end{cases}
	\end{align}
	Under the null hypothesis, we set $k^\ast=Tn$, so that only the first regime is observed.

	\begin{assumption}\label{assump1}
		Let $Q$ be the long-run covariance operator of
		$\{\widetilde Y_t^{(1)}\}_{t\in\mathbb Z}$, defined as in
		Equation~(\ref{eq11}). Assume that the jointly stationary process
		$\{(\widetilde Y_t^{(1)},\widetilde Y_t^{(2)})\}_{t\in\mathbb Z}$
		is $\rho$-mixing and that the following conditions hold for some
		$\delta>0$:
		\begin{enumerate}
			\item\label{assump1_p1}
			
			$\left\langle Q\mu_1,\mu_1\right\rangle_\bH>0.$

			\item\label{assump1_p2}
			
			$\sum_{m=1}^{\infty}
			[\alpha_{1,1}(m)]^{\delta/(2+\delta)}\le C,$
			where the mixing coefficients are those of
			$\{(\widetilde Y_t^{(1)},\widetilde Y_t^{(2)})\}_{t\in\mathbb Z}$ and the constant $C \in (0,\infty)$ is independent of $n$.
		\end{enumerate}
	\end{assumption}
	Part 1 of the Assumption is a simple non-degeneracy condition for the limiting Gaussian process. Part 2 requires a standard mixing assumption. The decay of the mixing coefficients is allowed to be fairly slow, because the data $\wt Y_t$ are bounded. Also, notice that strong mixing of $\{\wt Y_t\}_t$ is implied by strong mixing of the original data $\{ X_t\}_t$, because mixing coefficients are non-increasing under measurable transformations. 
	Now, under Assumption \ref{assump1} and the null hypothesis, Theorem \ref{th1} implies the FCLT
	\begin{align*}
		\left\{\frac{1}{\sqrt{n}}\sum_{t=1}^{\lfloor nr\rfloor}
		\big(\wt Y_t-\mu_t\big)\right\}_{r\in[0,T]}
		\rightsquigarrow
		\{B_Q(r)\}_{r\in[0,T]}
		\quad\text{in }D_\bH[0,T],
	\end{align*}
	where $\{B_Q(r)\}_{r\in[0,T]}$ is a $\bH$-valued Brownian motion with covariance operator $Q$. The following theorem establishes the limiting null distribution of $G_n$ under $H_0$; its proof is given in the supplement.
	
	\begin{theorem}\label{the_null}
		Suppose Assumption \ref{assump1} holds. Under $H_0$, we have
		\begin{align}
			G_{n}
			\stackrel{d}{\to}&  \sup_{\tau \in [1,T]}\sup_{r\in[1,\tau]}\frac{|B(r)-\frac{r}{\tau}B(\tau)|}{\sqrt{\int_{0}^{1}(B(s)-sB(1))^2ds}} \coloneqq \sup_{\tau \in [1,T]} \sup_{r \in [1,\tau]} G(\tau,r) \coloneqq G ,
		\end{align}
		where $\{B(t)\}_{t\in[0,T]}$ is an $\R$-valued standard Brownian motion.
	\end{theorem}
	
	\begin{remark}\label{remark_Wn}
		Under $H_0$, write $\mu=\mu_1$. The self-normalizer $W_n$ satisfies
		\begin{align}
			W_n\stackrel{d}{\to}
			\int_{0}^{1}\big\{\li\mu,B_Q(s)\ri_\bH-s\li\mu,B_Q(1)\ri_\bH\big\}^2ds
			=\li Q\mu,\mu\ri_\bH\int_{0}^{1}\{B(s)-sB(1)\}^2ds.
		\end{align}
		This means the distribution of $\sqrt{W_n}$ is asymptotically pivotal except for the long-run standard deviation $\sqrt{\li Q\mu,\mu\ri_\bH}$ which then cancels out the long-run standard deviation in the numerator, giving the fully pivotal limit $G$ of $G_n$.
	\end{remark}

	Next, we investigate the asymptotic properties of our monitoring statistic under the alternative. We introduce a sequence of local alternatives $\{H_{1,n}\}_{n=1}^\infty$ based on perturbing the mean of $\wt Y_t$ after the change. Specifically, consider a change at time $k^\ast=\lfloor nr_0\rfloor$ for some fixed $r_0\in(1,T)$ and 
	\begin{align}\label{eq_local_alternative}
		H_{1,n}:\quad
		\mu_t=\mu_1+n^{-\alpha}\Delta_\mu\id(t>k^\ast),
		\qquad t=1,\dots,Tn,
	\end{align}
	where $\Delta_\mu\in\bH\setminus\{0\}$ is fixed and $\alpha\in[0,\infty)$. Equivalently, the post-change mean embedding is $\mu_{Tn}=\mu_1+n^{-\alpha}\Delta_\mu$. The case $\alpha=0$ corresponds to a fixed alternative and we denote the time series underlying the post change segment as $\{\wt Y_t^A\}_{t\in\Z}$.  Theory for local alternatives can be mathematically derived under fairly general conditions on the pre- and post-change processes (using concepts such as uniform mixing coefficients). However, for the sake of transparency, we here restrict our attention to a more specialized, quite natural class of coupled alternatives. To understand its motivation, notice that for $\alpha>0$, by Lemma 4 in the Appendix of \cite{zhang2026doubly}, we have that $\wt Y_t^{(2)}$ has the same distribution as $(1-\delta_{n,t})\wt Y_t^{(1)}+\delta_{n,t}\wt Y_t^A$, where $\{\delta_{n,t}\}_{t\in\mathbb Z}$ is an i.i.d. sequence of $\operatorname{Bernoulli}(n^{-\alpha})$ random variables independent of the entire process $\{(\wt Y_t^{(1)},\wt Y_t^A)\}_{t\in\mathbb Z}$. For our alternatives, we now assume that $\wt Y_t^{(2)}=(1-\delta_{n,t})\wt Y_t^{(1)}+\delta_{n,t}\wt Y_t^A$ (not just in distribution, but on the same probability space), yielding a mixture of $\{\wt Y_t^{(1)}\}$ and $\{\wt Y_t^A\}$. Throughout this sequence, $H_{1,n}$ denotes the local alternative in
	Equation~(\ref{eq_local_alternative}) together with this mixture construction, and Assumption~\ref{assump1} is understood for the jointly stationary process $\{(\widetilde Y_t^{(1)},\widetilde Y_t^A)\}_{t\in\mathbb Z}$.

	To state weak limits, we also define $\eta:\Delta_T\to\R$ and 
	\begin{align}\label{eq_noncentral_limit}
		\mathcal G(\eta)
		=\sup_{\tau\in[1,T]}\sup_{r\in[1,\tau]}
		\frac{\left|B(r)-\frac{r}{\tau}B(\tau)+\eta(\tau,r)\right|}
		{\left[\int_0^1\{B(u)-uB(1)\}^2du\right]^{1/2}},
	\end{align}
	where $B$ is a standard Brownian motion. In particular, $\mathcal G(0)=G$. We also define for $(\tau,r)\in\Delta_T$
\[
\psi_{r_0}(\tau,r)
=[r-r_0]_+-\frac{r}{\tau}[\tau-r_0]_+,
\qquad [x]_+=\max\{x,0\}
\]
and (with $\theta_\Delta=\li\mu_1,\Delta_\mu\ri_\bH$)
\begin{align*}
\eta_{1/4}(\tau,r)
&=\frac{\|\Delta_\mu\|_\bH^2}
{\li Q\mu_1,\mu_1\ri_\bH^{1/2}}
\frac{[\tau-r_0]_+}{\tau}\psi_{r_0}(\tau,r),\\
\eta_{1/2}(\tau,r)
&=\frac{\theta_\Delta}
{\li Q\mu_1,\mu_1\ri_\bH^{1/2}}
\psi_{r_0}(\tau,r).
\end{align*}
	
	The following theorem characterizes the asymptotic power of $G_n$ under the fixed and local alternatives in Equation (\ref{eq_local_alternative}) and is proved in the supplement.
	
	\begin{theorem}\label{the_power}
		Suppose Assumption \ref{assump1} holds. Under $H_{1,n}$,
		\begin{enumerate}
			\item[(a)] If $\alpha=0$, then $G_n\stackrel{p}{\longrightarrow}\infty$.
			
			\item[(b)] For $\alpha>0$, define $\theta_\Delta=\li\mu_1,\Delta_\mu\ri_\bH$.
			\begin{enumerate}
				\item[(i)] If $\theta_\Delta=0$, $G_n\stackrel{p}{\longrightarrow}\infty$ if 
				$0<\alpha<\frac14$; $G_n\stackrel{d}{\longrightarrow}\mathcal G(\eta_{1/4})$ if $\alpha=\frac14$; and $G_n\stackrel{d}{\longrightarrow}G$ if $\alpha>\frac14$.
				
				\item[(ii)] If $\theta_\Delta\neq0$,  $G_n\stackrel{p}{\longrightarrow}\infty$ if 
				$0<\alpha<\frac12$; $G_n\stackrel{d}{\longrightarrow}\mathcal G(\eta_{1/2})$ if $\alpha=\frac12$; and $G_n\stackrel{d}{\longrightarrow}G$ if $\alpha>\frac12$.
			\end{enumerate}
		\end{enumerate}
		
	\end{theorem}
	
	Theorem \ref{the_power} identifies two local detection boundaries. When $\theta_\Delta\neq0$, the first-order signal yields the $n^{-1/2}$ boundary; when $\theta_\Delta=0$, the $a\|\Delta_\mu\|_\bH^2$ term in $h(a)$ from Equation (\ref{eq_power_signal}) becomes leading and yields the $n^{-1/4}$ boundary. Thus, the same monitoring-time variation of the projection that prevents $h(a)$ from vanishing identically and gives omnibus consistency under fixed alternatives also retains local power when the first-order alignment vanishes. In Section \ref{sec_power_enhancement} below, we will also see that the $n^{-1/4}$ boundary can be further improved by adding additional power-enhancement terms. We conclude with a small remark on weighted versions of our test statistic.

	\begin{remark} \label{rem:w}
		In addition to size and power, detection delay is an important criterion
		in sequential monitoring. \citet{kutta:dornemann:2025} showed that
		assigning extra weight to candidate splits close to the current monitoring time can
		shorten the detection delay of double-scan statistics such as ours. We briefly discuss how such a weighted version of our statistic would look like. For this purpose, define
		$T_{k,n}^w(s)=T_{k,n}(s)w((k-s)/n)$, where
		$w:(0,T-1]\to(0,\infty)$ is non-increasing, and set
		\begin{align*}
			G_n^w(k)
			&=
			\max_{s=n+1,\dots,k-1}
			\frac{|T_{k,n}^w(s)|}{\sqrt{W_n}}, \qquad
			G_n^w
			=
			\max_{n+2\leq k\leq Tn}G_n^w(k).
		\end{align*}
		Under appropriate additional conditions, $	G_n^w
		\stackrel{d}{\longrightarrow}
		\sup_{1<\tau\leq T}\sup_{1\leq r<\tau}
		G(\tau,r)w(\tau-r),$
		which remains pivotal and hence provides critical values for the
		weighted detector. The weights considered by
		\citet{kutta:dornemann:2025} are essentially of the form
		$w(x)=x^{-\beta}$, $\beta\in[0,1/2)$. The choice $\beta=0$ gives the
		unweighted statistic, whereas larger values shorten detection delays. The price of this approach is slightly lower power and a requirement of weaker temporal
		dependence, which would mean in this work faster-decaying mixing coefficients.
	\end{remark}

	\subsection{Recursive computation}\label{sec_iterative}
	
	The projection direction in Equation (\ref{eq_noisepro}) changes with the monitoring time. Nevertheless, the monitoring statistic can be updated without recomputing all kernel sums. Let
	\begin{align}\label{eq_iterative_scores}
		K_{ij}
		&=\li\wt Y_i,\wt Y_j\ri_\bH
		=f\big(d(X_i,X_j)\big),\quad
		Z_t^{(k)}
	=\left\langle\frac{U_k}{k},\wt Y_t\right\rangle_\bH
		=\frac{1}{k}\sum_{i=1}^{k}K_{it},
		\quad
		A_s^{(k)}=\sum_{t=1}^{s}Z_t^{(k)}.
	\end{align}
	Then Equation (\ref{eq_tkn}) can be written as
	\begin{align}\label{eq_iterative_tkn}
		T_{k,n}(s)
		=\frac{1}{\sqrt n}
		\left\{A_s^{(k)}-\frac{s}{k}A_k^{(k)}\right\}.
	\end{align}
	The historical scores also give
	\begin{align}\label{eq_iterative_wn}
		W_n
		=\frac{1}{n^2}\sum_{t=1}^{n}
		\left\{A_t^{(n)}-\frac{t}{n}A_n^{(n)}\right\}^2,
	\end{align}
	so $W_n$ is computed once from the historical sample and remains fixed throughout monitoring.
	
	\paragraph{Initialization.}
	For the historical sample, compute $h_t=\sum_{i=1}^{n}K_{it}$ and $Z_t^{(n)}=h_t/n$, $t=1,\dots,n$, whose cumulative sums give $W_n$ through Equation (\ref{eq_iterative_wn}). At the initial monitoring time $k_0=n+2$, compute $Z_t^{(k_0)}=k_0^{-1}\sum_{i=1}^{k_0}K_{it}$, $t=1,\dots,k_0$, and obtain $A_s^{(k_0)}$ by cumulative summation.
	
	\paragraph{Update.}
	Suppose $\{Z_t^{(k)}:1\leq t\leq k\}$ is available. When $X_{k+1}$ arrives, compute $K_{k+1,t}$ for $t=1,\dots,k+1$ and update
	\begin{align}\label{eq_iterative_update}
		Z_t^{(k+1)}
		&=\frac{k}{k+1}Z_t^{(k)}
		+\frac{1}{k+1}K_{k+1,t},
		\qquad t=1,\dots,k,\nonumber\\
		Z_{k+1}^{(k+1)}
		&=\frac{1}{k+1}\sum_{i=1}^{k+1}K_{i,k+1}.
	\end{align}
	A cumulative sum of the updated scores yields $A_s^{(k+1)}$ for every
	$s$. Hence,
	\begin{align}\label{eq_iterative_gn}
		G_n(k+1)
		=
		\max_{s=n+1,\dots,k}
		\frac{\left|A_s^{(k+1)}
			-\frac{s}{k+1}A_{k+1}^{(k+1)}\right|}
		{\sqrt{nW_n}}.
	\end{align}
	Each update requires $O(k)$ kernel evaluations and arithmetic
	operations. Therefore, for fixed $T$, computing the monitoring path
	through time $Tn$ requires $O(n^2)$ operations and $O(n)$ working
	storage, apart from storing the observations, whereas recomputing the
	projection scores directly at every monitoring time requires $O(n^3)$
	operations.
	
	\subsection{Power enhancement}
	\label{sec_power_enhancement}
	
	While Theorem~\ref{the_power} establishes omnibus consistency against general
	fixed alternatives, it also shows reduced local power when the change direction
	is orthogonal to the pre-change mean, that is, when
	$\theta_\Delta=\li\mu_1,\Delta_\mu\ri_\bH=0$. In this case, our
	statistic has an $n^{-1/4}$ rather than an optimal $n^{-1/2}$ detection boundary. Although exact orthogonality may be rare in practice, this result also
	suggests reduced finite-sample power when the angle between $\mu_1$ and
	$\Delta_\mu$ is close to $90^\circ$; we investigate this in our numerical experiments. To alleviate this
	problem, we propose a simple direction-agnostic power enhancement that recovers
	the $n^{-1/2}$ boundary up to logarithmic factors.  In the spirit of the power-enhancement
	principle of \citet{fan2015power}, we directly add a nonnegative,
	direction-agnostic detector $R_{k,n}(s)$ to the original numerator. Specifically, for a
	deterministic sequence $\lambda_n>0$, we define the power enhanced monitoring statistic as
	\begin{equation}
		G_n^{\mathrm{PE}}(k)
		=\max_{n+1\leq s<k}
		\frac{|T_{k,n}(s)|+\lambda_nR_{k,n}(s)}{\sqrt{W_n}},
		\qquad
		G_n^{\mathrm{PE}}=\max_{n+2\leq k\leq Tn}G_n^{\mathrm{PE}}(k).
		\label{eq_power_enhanced_statistic}
	\end{equation}
	Since $R_{k,n}(s)\geq0$, the enhanced statistic is pathwise at least as large
	as $G_n$. As shown in Theorem \ref{the_power_enhancement} below, the detector is asymptotically negligible under the null and diverges under alternatives.
	
	A natural candidate for this detector is the squared norm of the
	Hilbert-space CUSUM, $\|C_{k,n}(s)\|_\bH^2$.  Directly using this
	quantity, however, may create a finite-sample distortion due to the same reason as stated in Remark \ref{rmk_2}. We therefore use the
	adjusted squared CUSUM
	\begin{equation}
		R_{k,n}(s)
		:=\left[
		\|C_{k,n}(s)\|_\bH^2
		-f(0)\frac{s(k-s)}{nk}
		\right]_+,
		\qquad [x]_+=\max\{x,0\}.
		\label{eq_pe_component}
	\end{equation}
	The following theorem shows the asymptotic properties of $G_n^{\mathrm{PE}}(k)$ and is proved in the supplement.
	\begin{theorem}
		\label{the_power_enhancement}
		Suppose Assumption~\ref{assump1} holds and $\lambda_n>0$ satisfies
		$\lambda_n\to0$.
		\begin{enumerate}
			\item[(a)] Under the null hypothesis: $
			\max_{n+2\leq k\leq Tn}
			\big\{G_n^{\mathrm{PE}}(k)-G_n(k)\big\}=o_p(1),$ and $G_n^{\mathrm{PE}}\stackrel{d}{\longrightarrow}G.$
			
			\item[(b)] Under the alternative in Equation~(\ref{eq_local_alternative}) with $n^{-\alpha}$ replaced by a
			sequence $c_n\downarrow0$.
			\begin{enumerate}
				\item[(i)] If $\theta_\Delta=0$, $\sqrt n\,\lambda_n\to\infty$
				and $n\lambda_nc_n^2\to\lambda\in[0,\infty)$, then
				\[
				G_n^{\mathrm{PE}}\stackrel{d}{\longrightarrow}
				G^{\mathrm{PE}}(\lambda),
				\]
				where, for $\lambda\geq0$,
\[
G^{\mathrm{PE}}(\lambda)
=\sup_{\tau\in[1,T]}\sup_{r\in[1,\tau]}
\frac{
\left|B(r)-\frac{r}{\tau}B(\tau)\right|
+\lambda\{\li Q\mu_1,\mu_1\ri_\bH\}^{-1/2}
\|\Delta_\mu\|_\bH^2\psi_{r_0}(\tau,r)^2
}{
\left[\int_0^1\{B(u)-uB(1)\}^2du\right]^{1/2}
}.
\]
Here, $B$ is a standard Brownian motion on $[0,T]$, and
$\psi_{r_0}$ is defined before Theorem~\ref{the_power}.
				
				\item[(ii)] Without any restriction on $\theta_\Delta$, if
				$n\lambda_nc_n^2\to\infty$, then
				$G_n^{\mathrm{PE}}\stackrel{p}{\longrightarrow}\infty$.
			\end{enumerate}
		\end{enumerate}
	\end{theorem}
	
	If $\theta_\Delta\neq0$, the
	sharper $n^{-1/2}$ detection boundary follows directly from
	Theorem~\ref{the_power}, because $G_n^{\mathrm{PE}}\geq G_n$. If
	$\theta_\Delta=0$, the condition
	$n\lambda_nc_n^2\to\infty$ allows $c_n$ to decay at a close-to-root-$n$ rate. For
	example, with $\lambda_n=(\log n)^{-\kappa}$, this condition holds whenever $c_n\gg\frac{(\log n)^{\kappa/2}}{\sqrt n}$. Thus, the detection
	boundary for orthogonal changes improves from $n^{-1/4}$ to nearly $n^{-1/2}$ with the power enhanced monitoring statistic $G_n^{\mathrm{PE}}(k)$. In the simulations in Section~\ref{sec3simu}, we use
	$\lambda_n=(\log n)^{-\kappa}$ with $\kappa=3$.

	\section{Open-End Monitoring}\label{sec_open}
	
	In this section, we consider open-end monitoring, under which the procedure continues indefinitely if no change is detected. Open-end monitoring is a standard problem-formulation already since \cite{chu1996monitoring} and we refer to \citet{chan2021self},  \citet{gossmann2021new} and \cite{aue:kirch:2024} for examples. Open-end monitoring can be viewed as the monitoring setup in Section~\ref{sec2_2} with $T=\infty$. 
	This setting is useful, because the analyst does not have to commit to a fixed endpoint of monitoring, before the procedure is even launched. In the following,
	we retain the detector $G_n(k)$ and the historical self-normalizer $W_n$, but introduce a new monitoring-time-dependent weight function $g$. To avoid confusion with the weights in Remark \ref{rem:w}, we call $g$ the ``(temporal) discounting". 
	Now, define
	\begin{align}\label{eq_open_stat}
		G_{n,g}
		=\sup_{k\geq n+2}\frac{1}{g(k/n)}G_n(k)
		=\sup_{k\geq n+2}\frac{1}{g(k/n)}
		\max_{s=n+1,\ldots,k-1}\frac{|T_{k,n}(s)|}{\sqrt{W_n}}.
	\end{align}
	As in \citet{gossmann2021new}, temporal discounting is needed to obtain a nondegenerate open-end limiting null distribution. Throughout this section, we assume that the function $g$  satisfies the following assumption.
	
	\begin{assumption}\label{assump_open_weight}
		The function $g:[1,\infty)\to(0,\infty)$ is continuous and nondecreasing, and
		\begin{align}
			\lim_{\tau\to \infty}	\frac{\sqrt{\tau\log\log\tau}}{g(\tau)}& = 0
			, \label{eq_g1}\\
			\lim_{T\to\infty}\sum_{\ell=0}^{\infty}
			\frac{2^\ell T}{g(2^\ell T)^2}&=0. \label{eq_g2}
		\end{align}
	\end{assumption}
	
	\begin{remark}\label{remark_open_weight}
		Assumption ~\ref{assump_open_weight} is satisfied by the linear discounting $g(\tau)=\tau$  used by \citet{chan2021self}. It is also satisfied by the function $g_\gamma$ considered in \citet{gossmann2021new} which is $g_\gamma(\tau)=\tau\max\{((\tau-1)/\tau)^\gamma,\epsilon\}$, for $\tau\geq1$, $0\leq\gamma<1/2$ and $\epsilon>0$ . 
	\end{remark}
	
	Let $B(\cdot)$ be a standard Brownian motion on $[0,\infty)$ and define
	$\Delta_\infty=\{(\tau,r):1\leq r\leq\tau<\infty\}$. For any bounded function $\eta:\Delta_\infty\to\R$, define
	\begin{align}\label{eq_open_limit_eta}
		\mathcal G_g(\eta)
		=\frac{1}{\left[\int_0^1\{B(u)-uB(1)\}^2du\right]^{1/2}}
		\sup_{\tau\geq1}\sup_{1\leq r\leq\tau}
		\frac{\left|B(r)-\frac r\tau B(\tau)+\eta(\tau,r)\right|}{g(\tau)},
	\end{align}
	and write $G_g=\mathcal G_g(0)$.
	The following theorem establishes the limiting null distribution of $G_{n,g}$ under $H_0$; its proof is given in the supplementary material.
	
	\begin{theorem}\label{the_open_null}
		Suppose Assumptions~\ref{assump1} and \ref{assump_open_weight} hold. Under $H_0$ we have
		$G_{n,g}\stackrel{d}{\longrightarrow}G_g$.
		
	\end{theorem}

	Approximating the limiting quantiles of $G_g$ is in principle not difficult; one can approximate the Brownian motion by a partial sum process and replace maximization over the infinite domain $\Delta_\infty$ by $\Delta_T$  for some sufficiently large  $T \in \mathbb{R}$. However, what $T$ is ``sufficiently large" in fact depends on the weight function  $g$, where slower growth requires larger $T$. To circumvent a choice of $T$ entirely, we provide a rescaled version of our limit to the unit interval in the next lemma. Its proof is given in the supplement.

	\begin{lemma}\label{lem_open_fixed_interval}
		Let $W$ and $W_0$ be independent standard Brownian motions on $[0,1]$. Then
		\begin{align}\label{eq_open_fixed_interval}
			G_g\stackrel{d}{=}
			\frac{1}{\left[\int_0^1\{W_0(u)-uW_0(1)\}^2du\right]^{1/2}}
			\sup_{0<t\leq v\leq1}
			\frac{|W(v)-W(t)|}{v g(1/t)}.
		\end{align}
	\end{lemma}
	
	\par\medskip
	For the asymptotic results under alternative, we use the alternatives $H_{1,n}$ from Equation
	(\ref{eq_local_alternative}), with a fixed $r_0>1$.
	
	\begin{theorem}\label{the_open_power}
		Suppose the conditions of Theorem~\ref{the_open_null} hold. Under $H_{1,n}$,
		\begin{enumerate}
			\item[(a)] If $\alpha=0$, then $G_{n,g}\stackrel{p}{\longrightarrow}\infty$.
			\item[(b)] If $\alpha>0$, define $\theta_\Delta=\langle\mu_1,\Delta_\mu\rangle_\bH$.
			\begin{enumerate}
				\item[(i)] If $\theta_\Delta=0$, then $G_{n,g}\stackrel{p}{\longrightarrow}\infty$ if $0<\alpha<1/4$; $G_{n,g}\stackrel{d}{\longrightarrow}\mathcal G_g(\eta_{1/4})$ if $\alpha=1/4$; and $G_{n,g}\stackrel{d}{\longrightarrow}G_g$ if $\alpha>1/4$.
				\item[(ii)] If $\theta_\Delta\neq0$, then $G_{n,g}\stackrel{p}{\longrightarrow}\infty$ if $0<\alpha<1/2$; $G_{n,g}\stackrel{d}{\longrightarrow}\mathcal G_g(\eta_{1/2})$ if $\alpha=1/2$; and $G_{n,g}\stackrel{d}{\longrightarrow}G_g$ if $\alpha>1/2$.
			\end{enumerate}
		\end{enumerate}
		Here, $\eta_{1/4}$ and $\eta_{1/2}$ are given by the same formulas as before Theorem~\ref{the_power}, now for $(\tau,r)\in\Delta_\infty$.
	\end{theorem}
	
	Thus, for a fixed $r_0>1$, the open-end procedure retains the same local detection boundaries as the closed-end procedure, depending on whether the change direction $\Delta_\mu$ is orthogonal to the pre-change mean embedding $\mu_1$. The proof is given in the supplementary material.
	
	\paragraph{Recursive implementation.}
	The open-end procedure uses exactly the recursive quantities in Section~\ref{sec_iterative}. At time $k$, compute $G_n(k)$ from Equations~(\ref{eq_iterative_update})--(\ref{eq_iterative_gn}), divide by $g(k/n)$, and stop when this value exceeds the open-end critical value. No additional kernel sums are needed. Although each update still costs $O(k)$ operations, an open-end procedure has no deterministic finite total cost because its stopping horizon is unbounded.
	
	\section{Simulation results}\label{sec3simu}
		
		We investigate the finite-sample performance of the proposed
		procedures in two simulation studies, one for distribution-valued time series in this section and a second one for graph-valued time series in the supplementary material.
		Throughout, we consider a monitoring horizon of $T=5$, historical sample sizes
		$n\in\{300,500\}$, and the nominal level $\alpha=0.05$. Empirical null
		rejection probabilities are based on $3{,}000$ Monte Carlo replications,
		whereas power is estimated from $500$ replications under each alternative.

		\paragraph{Benchmarks.}
		To the best of our knowledge there currently does not exist another method 
		for sequential change detection in the distribution of dependent object-valued time series.
		We therefore establish a reasonable benchmark as follows: We combine our double kernel embedding with an established procedure for monitoring mean changes in a
		Hilbert space. Recall that the double embedding maps
		each original observation $X_i$ to an element $\widetilde Y_i$ of an RKHS,
		thereby reducing the original problem of detecting a distributional change in $(X_i)_{i\in\mathbb N}$ to detecting a mean change in
		$(\widetilde Y_i)_{i\in\mathbb N}$. Therefore, the sequential CUSUM procedure of
		\citet{kutta2025monitoring}, henceforth referred to as KK, can be applied
		directly to the twice embedded observations. We compare our self-normalized
		detector (SN) and its power-enhanced version (SN-PE) with KK in terms of null
		rejection probabilities and power. SN-PE is described in Section \ref{sec_power_enhancement}, where we have chosen $\lambda_n = (\log(n))^{-3}$. In addition, we also compare with the three	monitoring procedures of \citet{boniece2026sequential}, denoted by $D_1$, $D_2$,
		and $D_3$, respectively. These procedures are designed for independent random objects and are able to detect distributional changes using
		degenerate two-sample $U$-statistics. 
		
		\subsection{Distribution-valued observations}\label{sec3simu_distribution}
		
		\paragraph{Data generation.}
		Our setting is inspired by \citet{jiang2024two} and
		\citet{zhang2025change}, who consider retrospective tests for
		non-Euclidean time series. Let $q\geq 0$ and consider
		\begin{equation}
			U_t^{(q)}
			=\frac{1}{\sqrt{q+1}}\sum_{j=0}^{q}\varepsilon_{t-j},
			\qquad
			\varepsilon_t\stackrel{\mathrm{iid}}{\sim}\mathcal N(0,1).
			\label{eq:ma-process}
		\end{equation}
		Let $\beta_t\stackrel{\mathrm{iid}}{\sim}\operatorname{Bernoulli}(\delta)$ and
		$E_t\stackrel{\mathrm{iid}}{\sim}\operatorname{Exp}(1)$, independently of
		each other and of the innovations in \eqref{eq:ma-process}. For
		$k^*=\lfloor r_0n\rfloor$, define
		\begin{equation}
			X_t=\mathcal N(m_t,1),
			\qquad
			m_t=
			\begin{cases}
				uU_t^{(q)}, & 1\leq t\leq k^*,\\
				(1-\beta_t)uU_t^{(q)}+\beta_tE_t, & k^*<t\leq 5n.
			\end{cases}
			\label{eq:distribution-dgp}
		\end{equation}
		The null hypothesis corresponds to $\delta=0$, and larger values of $\delta$ lead deeper into the alternative. Change locations are
		$r_0\in\{1.5,3,4.5\}$, representing early, intermediate, and late changes,
		respectively. Under the alternative, each post-change observation is
		independently replaced by $\mathcal N(E_t,1)$ with probability $\delta$.
		The observation space is equipped with the $2$-Wasserstein metric. Since all
		normal distributions in \eqref{eq:distribution-dgp} have variance one,
		$d(X_i,X_j)=|m_i-m_j|$. For the second embedding, we use the kernel
		$\mathcal K(Y_i,Y_j)=\exp\{-d(X_i,X_j)\}$.\\
		The model contains two parameters of particular interest. First, $q$
		determines the strength of temporal dependence, with larger values
		corresponding to stronger dependence. Its effect is particularly relevant
		under the null hypothesis, where it may affect the nominal approximation.  We consider $q \in \{1,3,6\}$ below.
		Second, $u>0$ determines the distributions of the embedded observations
		before and after the change. In particular, it influences the angle between
		the pre-change mean $\mu_1$ and the change direction $\Delta_\mu$. Angles
		closer to $90^\circ$ are especially challenging for SN. We consider
		$u\in\{0.12,0.44,0.65,0.90\}$, corresponding approximately to angles of
		$50^\circ$, $60^\circ$, $70^\circ$, and $80^\circ$. Since $u$ also affects
		the magnitude of the change, the power comparisons below are indexed by
		$\|\Delta_\mu\|$ rather than by $\delta$.
		
		\paragraph{Parameter choices and quantile calibration.}
		The KK detector involves the choice of two parameters, which we set to
		$\gamma=0.3$ and $\zeta=0.05$, following the simulation design of
		\citet{kutta2025monitoring}. Both SN and KK require approximations of their
		limiting quantiles. For SN, the limiting distribution is pivotal. The same
		critical value can therefore be used for SN and SN-PE across all considered
		values of $n$, $q$, and $u$. We approximate this quantile once by simulating
		the corresponding one-dimensional limit. Based on $20{,}000$
		replications and $T=5$, we obtain $q_{\mathrm{SN},5}=13.255$.\\
		For KK, the limiting distribution depends on the spectrum of the long-run
		covariance operator. We consider an oracle-type calibration, referred to as
		oracle KK below, and a practical calibration. For oracle KK, the spectrum is
		estimated once for each $(q,u)$ setting from an independent null pilot sample
		of length $1{,}500$, using the known dependence structure and the exact
		long-run covariance trace. The resulting Gaussian critical value is
		approximated from $5{,}000$ replications and then reused throughout the
		simulation. For practical KK, the long-run covariance operator is re-estimated
		from the training sample in every Monte Carlo replication, using the approach
		of \citet{rice:shang:2017} with a plug-in bandwidth and Bartlett weights. The
		corresponding critical value is based on $500$ conditional Gaussian bootstrap
		replications. This comparison separates the finite-sample effect of long-run
		covariance estimation from the intrinsic behavior of the KK statistic.\\
		For $D_1,D_2,D_3$, we use $h(X_i,X_j)=\exp\{-d(X_i,X_j)\}$ and set
		$\beta=0.5$, with $c_0=0.5$ and $\gamma=0.51$ for $D_3$.
		We use the calibration procedure in Section~4.4 of \citet{boniece2026sequential}.
		
		\paragraph{Nominal approximation.}
		Table~\ref{Tab:1} reports empirical rejection probabilities under the null
		hypothesis for the nominal level $5\%$. SN provides a reasonable
		nominal approximation across all considered values of $u$ and $q$, with a
		slight improvement as the training-sample size increases. The rejection
		probabilities of SN-PE are slightly higher than those of SN. This is expected,
		since the SN-PE statistic is pathwise no smaller than the SN statistic, while
		both procedures use the same limiting critical value. The nominal
		approximation of oracle KK is comparable to that of SN. Oracle KK is, however,
		slightly more liberal than SN in several settings, particularly for $n=500$
		and $q\in\{3,6\}$, although the differences remain small.\\
		By contrast, practical KK, which estimates the long-run covariance operator
		from the training sample, exhibits substantial level inflation. The distortion
		becomes more pronounced as $q$ increases. For $q=6$ and $n=300$, the empirical
		rejection probability is approximately $16\%$ across the different
		$u$-scenarios. Increasing the training-sample size to $n=500$ improves the
		approximation, but the rejection probability remains approximately $13\%$ for
		$q=6$, almost three times the nominal level. The inflation is less severe, but
		still considerable, for $q=3$, and smaller, though still visible, for $q=1$.
		These findings are consistent with evidence from the existing literature. For
		example, \citet[see their Tables~2 and~3]{gossmann2021new} demonstrate that,
		even for real-valued data, estimation of the long-run variance can
		substantially impair finite-sample size control relative to oracle
		calibration. The effect appears to be even more pronounced in the present
		Hilbert-space setting, highlighting an important practical advantage of
		self-normalization.
		
		The rejection probabilities of $D_1,D_2,D_3$ are substantially higher:
		they range from $21.4\%$ to $37.3\%$ for $q=1$, from $57.9\%$ to $89.9\%$ for
		$q=3$ and from $86.8\%$ to $99.8\%$ for $q=6$. In every setting, $D_3$ has
		the largest rejection probability, followed by $D_2$ and $D_1$.
		Increasing $n$ from $300$ to $500$ does not remedy these distortions.
		These results concern the robustness of their method outside the independent-data setting covered by their theory.
		
		\begin{table}[h]
			\centering
			\caption{Empirical rejection probabilities under the null hypothesis at the
				nominal level $\alpha=0.05$, based on $3{,}000$ Monte Carlo replications.}
			\label{Tab:1}
			\small
			\setlength{\tabcolsep}{5pt}
			\begin{tabular}{@{}cccccccccc@{}}
				\toprule
				$n$ & $q$ & $u$ & SN & SN-PE & KK (oracle) & KK (practical) & $D_1$ & $D_2$ & $D_3$ \\
				\midrule
				300 & 1 & 0.12 & 0.045 & 0.045 & 0.044 & 0.078 & 0.248 & 0.282 & 0.361 \\
				300 & 1 & 0.44 & 0.044 & 0.045 & 0.044 & 0.073 & 0.234 & 0.270 & 0.334 \\
				300 & 1 & 0.65 & 0.046 & 0.048 & 0.045 & 0.076 & 0.229 & 0.257 & 0.313 \\
				300 & 1 & 0.90 & 0.044 & 0.047 & 0.045 & 0.072 & 0.214 & 0.244 & 0.297 \\
				\addlinespace
				300 & 3 & 0.12 & 0.055 & 0.056 & 0.057 & 0.104 & 0.608 & 0.722 & 0.872 \\
				300 & 3 & 0.44 & 0.051 & 0.058 & 0.052 & 0.106 & 0.609 & 0.709 & 0.866 \\
				300 & 3 & 0.65 & 0.052 & 0.057 & 0.050 & 0.109 & 0.593 & 0.701 & 0.853 \\
				300 & 3 & 0.90 & 0.051 & 0.057 & 0.050 & 0.107 & 0.579 & 0.683 & 0.834 \\
				\addlinespace
				300 & 6 & 0.12 & 0.066 & 0.068 & 0.063 & 0.154 & 0.870 & 0.951 & 0.996 \\
				300 & 6 & 0.44 & 0.057 & 0.064 & 0.060 & 0.168 & 0.881 & 0.953 & 0.996 \\
				300 & 6 & 0.65 & 0.053 & 0.063 & 0.061 & 0.162 & 0.876 & 0.956 & 0.996 \\
				300 & 6 & 0.90 & 0.050 & 0.062 & 0.061 & 0.167 & 0.868 & 0.946 & 0.994 \\
				\midrule
				500 & 1 & 0.12 & 0.047 & 0.047 & 0.046 & 0.078 & 0.252 & 0.293 & 0.373 \\
				500 & 1 & 0.44 & 0.047 & 0.047 & 0.049 & 0.075 & 0.241 & 0.279 & 0.343 \\
				500 & 1 & 0.65 & 0.046 & 0.048 & 0.050 & 0.076 & 0.227 & 0.266 & 0.330 \\
				500 & 1 & 0.90 & 0.047 & 0.051 & 0.053 & 0.076 & 0.216 & 0.246 & 0.310 \\
				\addlinespace
				500 & 3 & 0.12 & 0.050 & 0.051 & 0.065 & 0.090 & 0.640 & 0.768 & 0.899 \\
				500 & 3 & 0.44 & 0.053 & 0.056 & 0.060 & 0.090 & 0.638 & 0.753 & 0.887 \\
				500 & 3 & 0.65 & 0.052 & 0.055 & 0.059 & 0.092 & 0.626 & 0.739 & 0.879 \\
				500 & 3 & 0.90 & 0.051 & 0.057 & 0.057 & 0.094 & 0.614 & 0.719 & 0.857 \\
				\addlinespace
				500 & 6 & 0.12 & 0.054 & 0.055 & 0.070 & 0.123 & 0.896 & 0.967 & 0.998 \\
				500 & 6 & 0.44 & 0.049 & 0.055 & 0.066 & 0.125 & 0.901 & 0.970 & 0.996 \\
				500 & 6 & 0.65 & 0.049 & 0.054 & 0.065 & 0.129 & 0.897 & 0.968 & 0.997 \\
				500 & 6 & 0.90 & 0.047 & 0.052 & 0.067 & 0.135 & 0.896 & 0.964 & 0.996 \\
				\bottomrule
			\end{tabular}
		\end{table}
		
		\paragraph{Power.}
		For the power comparison, we consider only SN, SN-PE, and oracle KK.
		Given the severe size distortions of $D_1,D_2,D_3$, we do not conduct power
		experiments for these procedures. We omit
		practical KK to separate the intrinsic power properties of the detectors from
		the finite-sample effects of long-run covariance estimation. In an additional
		pilot study for $n=300$ and $q\in\{1,3\}$, simulations (not reported
		here) indicate that the size-adjusted power of practical KK is nearly identical
		to that of oracle KK.  Within each of the following alternative scenarios,
		the procedures are evaluated on the same simulated paths.
		Figures~\ref{fig:power-n500-angle50}--\ref{fig:power-n500-angle80} report the
		results for $n=500$; the corresponding figures for $n=300$ are reported in
		the supplementary material. The figures show size-adjusted empirical power against the change
		magnitude $\|\Delta_\mu\|$. For each fixed value of $u$, the change magnitude
		is varied through $\delta$ in \eqref{eq:distribution-dgp}, without changing the
		angle between $\mu_1$ and $\Delta_\mu$. The results reveal several distinct,
		and partly overlapping, effects.\\
		First, KK performs best for early changes and worst for late changes. The main
		reason is that KK compares the historical mean with the mean of all collected
		monitoring observations. If the change occurs late, the monitoring sample
		contains many pre-change observations, which dilutes the signal. Moreover, the
		boundary weighting with $\gamma=0.3$ gives KK relatively greater sensitivity
		at early monitoring times. By contrast, SN scans jointly over the monitoring
		time and potential change locations and is therefore less sensitive to the
		change location. These different scan types have been discussed in detail in the review article by \citet{aue:kirch:2024} (they are referred to as ``CUSUM'' and ``full-CUSUM'' scans by these authors, respectively).
		SN's power still also decreases for very late changes
		because only a small post-change sample remains available.\\
		Second, the change in $u$, and hence in the associated angle between $\mu_1$
		and $\Delta_\mu$, has an important effect on the relative performance of the
		procedures. As the angle increases, the projection underlying SN removes a
		larger part of the signal. At an angle of approximately $80^\circ$, KK
		substantially outperforms standard SN for early and intermediate changes,
		whereas its advantage disappears for late changes. At $70^\circ$, neither
		procedure uniformly dominates: KK performs better for early changes, while SN
		performs somewhat better for intermediate changes and much better for late changes.
		At angles of $60^\circ$ and $50^\circ$, SN typically outperforms KK; results are comparable for early changes at $60^\circ$, but at $50^\circ$, SN dominates by a large margin in any scenario. A plausible explanation is
		that, at these smaller angles, projection reduces the noise more strongly than
		it reduces the signal. This finite-sample dimension-reduction effect is not
		directly captured by our asymptotic theory, but may be very relevant in
		applications. We did not investigate the effect of even smaller angles, even though we expect the advantage of SN there to be even larger.\\
		As expected, stronger dependence also tends to reduce power: across the panels, power for
		$q=3$ is generally lower than for $q=1$, particularly at larger angles and for
		intermediate or late changes. The qualitative rankings of the procedures,
		however, remain broadly similar.\\
		Finally, SN-PE produces moderate but systematic power gains over standard SN,
		particularly at angles of $70^\circ$ and $80^\circ$. The two procedures are
		nearly indistinguishable at $50^\circ$ and gains are small at
		$60^\circ$. 
		
		This corresponds to our theoretical expectation that for narrower angles, the SN part dominates the power behavior, while for large angles, the power enhancement makes a bigger contribution. The power enhancement can be substantial. For example, at $80^\circ$,
		$r_0=3$, $n=500$, and $q=1$, the maximum power difference exceeds
		$30$ percentage points.
	\begingroup
	\captionsetup{font={footnotesize,singlespacing},skip=3pt}
		
		\begin{figure}[p]
			\centering
			\includegraphics[width=.93\linewidth]
			{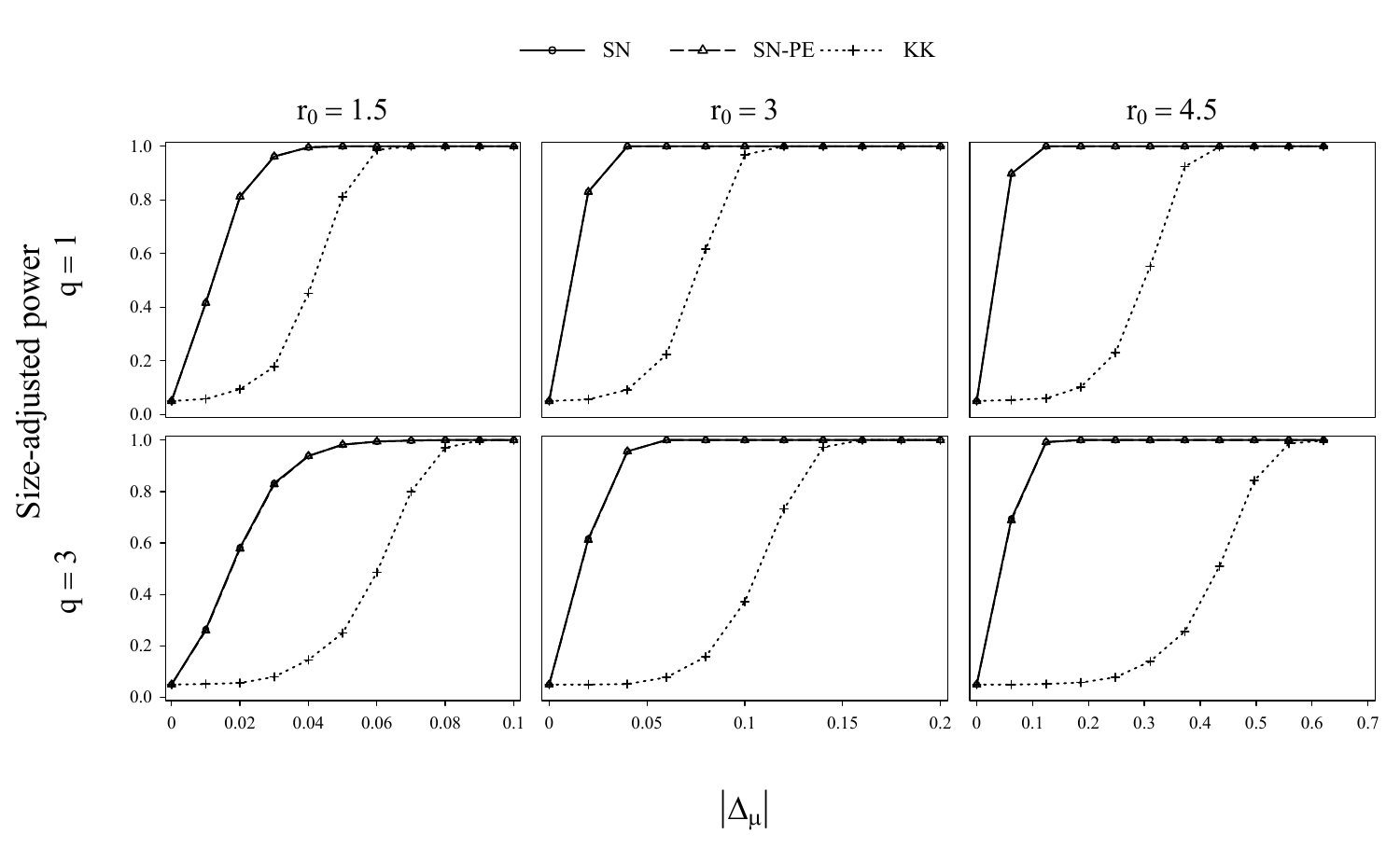}
			\caption{Size-adjusted empirical power for $n=500$ and $u=0.12$,
				corresponding to an angle of approximately $50^\circ$. Columns correspond
				to $r_0=1.5$, $3$, and $4.5$, while rows correspond to $q=1$ and $q=3$.
				Each panel compares SN, SN-PE, and oracle KK as a function of the change
				magnitude $\|\Delta_\mu\|$.}
			\label{fig:power-n500-angle50}
		\end{figure}

		\begin{figure}[p]
			\centering
			\includegraphics[width=.93\linewidth]
			{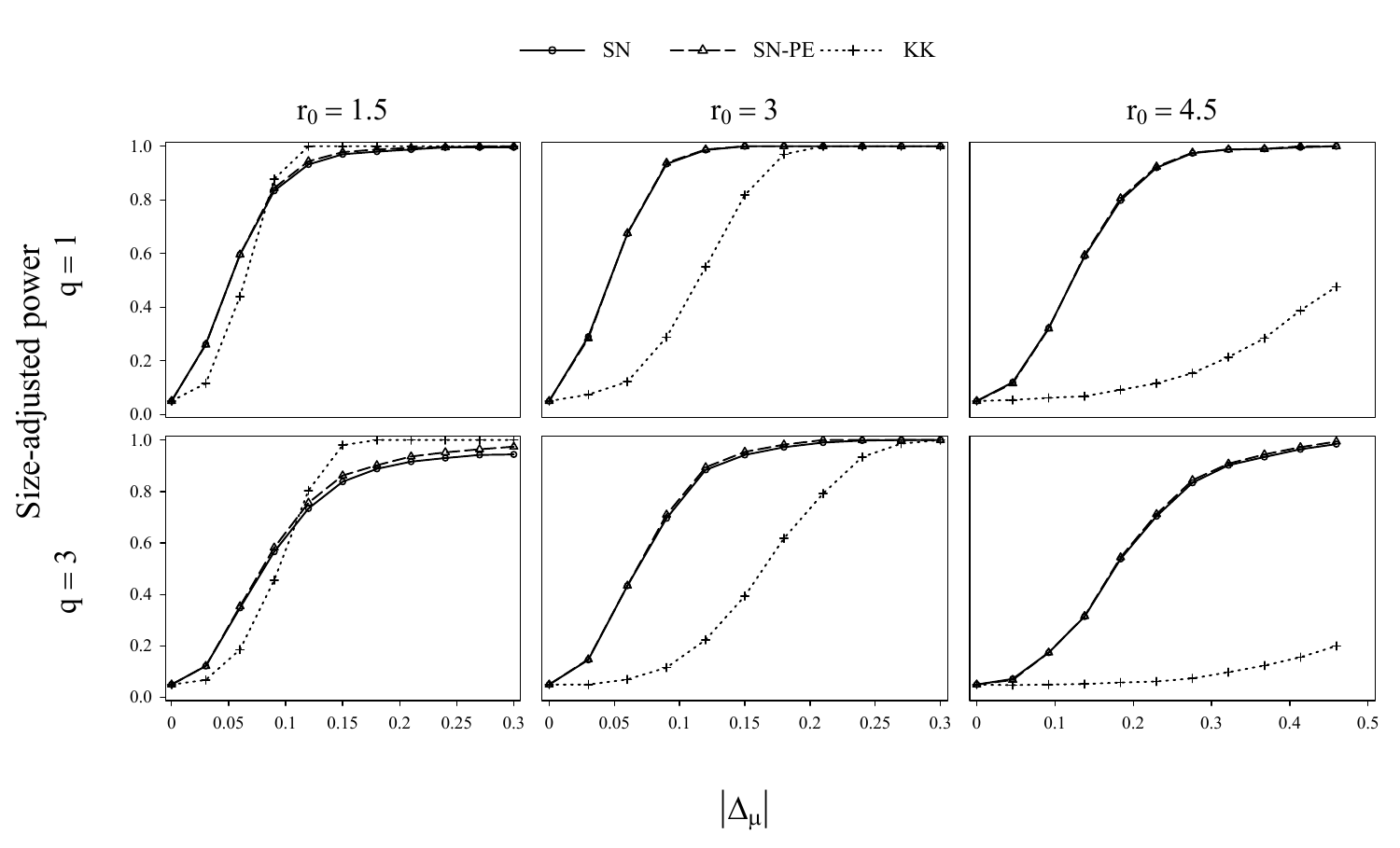}
			\caption{Size-adjusted empirical power for $n=500$ and $u=0.44$,
				corresponding to an angle of approximately $60^\circ$. Columns correspond
				to $r_0=1.5$, $3$, and $4.5$, while rows correspond to $q=1$ and $q=3$.
				Each panel compares SN, SN-PE, and oracle KK as a function of the change
				magnitude $\|\Delta_\mu\|$.}
			\label{fig:power-n500-angle60}
		\end{figure}
		
		\clearpage
		
		\begin{figure}[p]
			\centering
			\includegraphics[width=.93\linewidth]
			{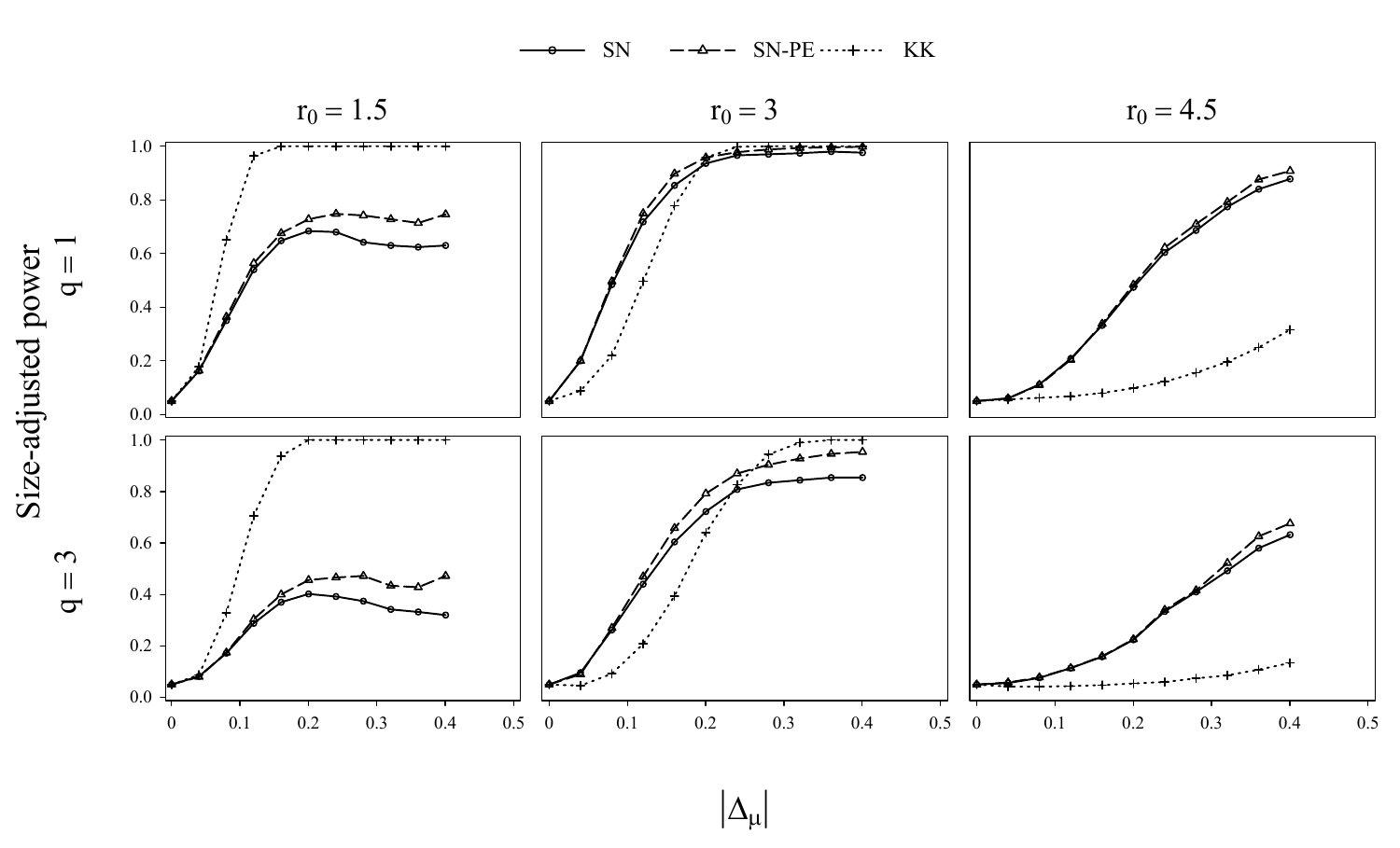}
			\caption{Size-adjusted empirical power for $n=500$ and $u=0.65$,
				corresponding to an angle of approximately $70^\circ$. Columns correspond
				to $r_0=1.5$, $3$, and $4.5$, while rows correspond to $q=1$ and $q=3$.
				Each panel compares SN, SN-PE, and oracle KK as a function of the change
				magnitude $\|\Delta_\mu\|$.}
			\label{fig:power-n500-angle70}
		\end{figure}

		\begin{figure}[p]
			\centering
			\includegraphics[width=.93\linewidth]
			{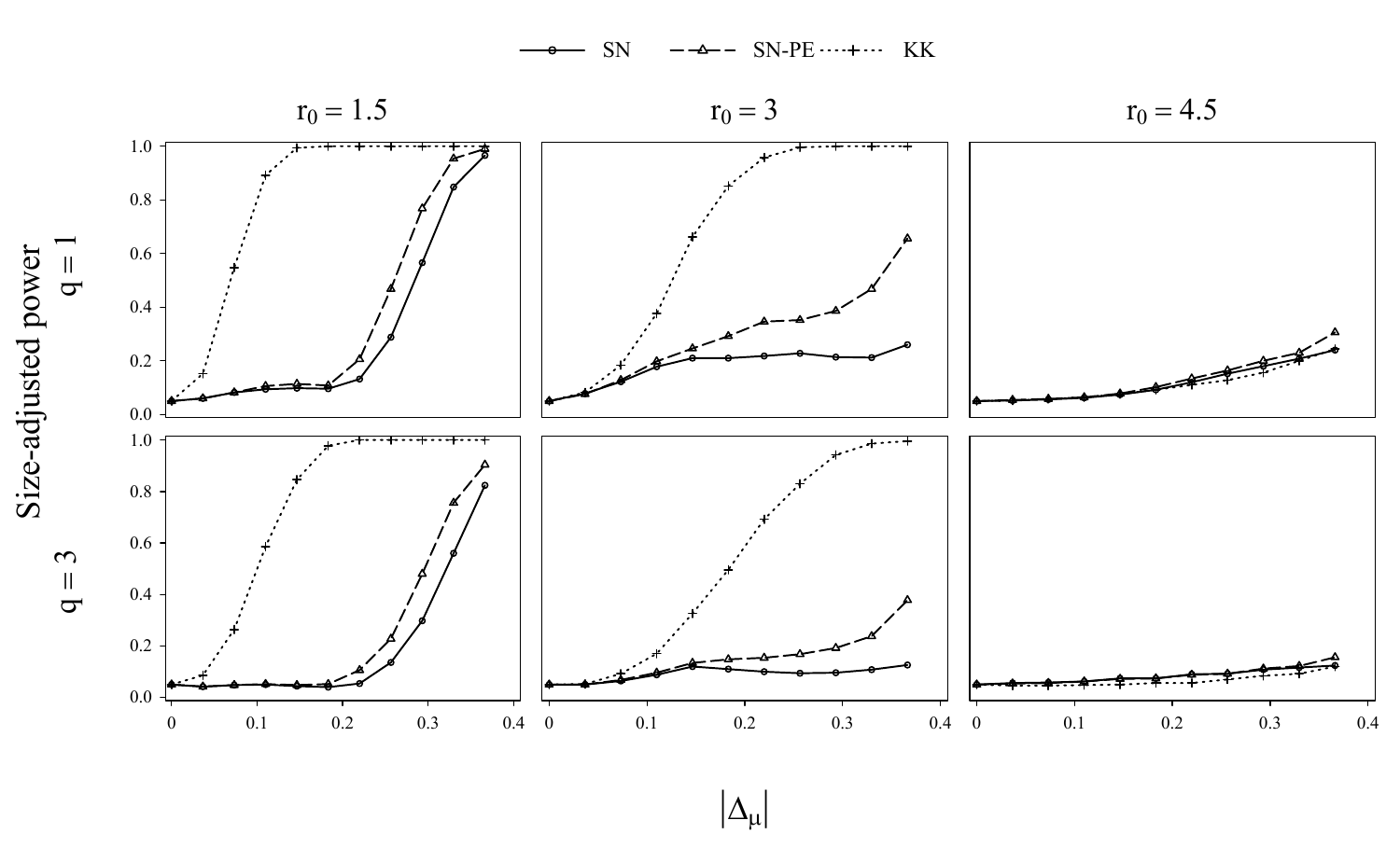}
			\caption{Size-adjusted empirical power for $n=500$ and $u=0.90$,
				corresponding to an angle of approximately $80^\circ$. Columns correspond
				to $r_0=1.5$, $3$, and $4.5$, while rows correspond to $q=1$ and $q=3$.
				Each panel compares SN, SN-PE, and oracle KK as a function of the change
				magnitude $\|\Delta_\mu\|$.}
			\label{fig:power-n500-angle80}
		\end{figure}

	\clearpage
	\endgroup

	\section{Application to Stock-Return Distributions}\label{sec4real}
	
	We apply the proposed monitoring procedure to monthly distributions of stock returns. Daily adjusted prices for stocks in the S\&P 500 were obtained from Yahoo Finance over January 2012--December 2025, and adjusted log returns were computed for each trading day. For each calendar month, we pool the available daily returns across stocks and regard their distribution as one observation $X_t$. This gives 168 monthly observations. We consider two monitoring frames, both with $T=2$. In the first, January 2012--December 2018 serves as the stable historical period ($n=84$), and January 2019--December 2025 is monitored. In the second, the stable historical period is January 2012--December 2016 ($n=60$), and January 2017--December 2021 is monitored.
	
\begin{figure}[H]
		\centering
		\includegraphics[width=0.92\linewidth]{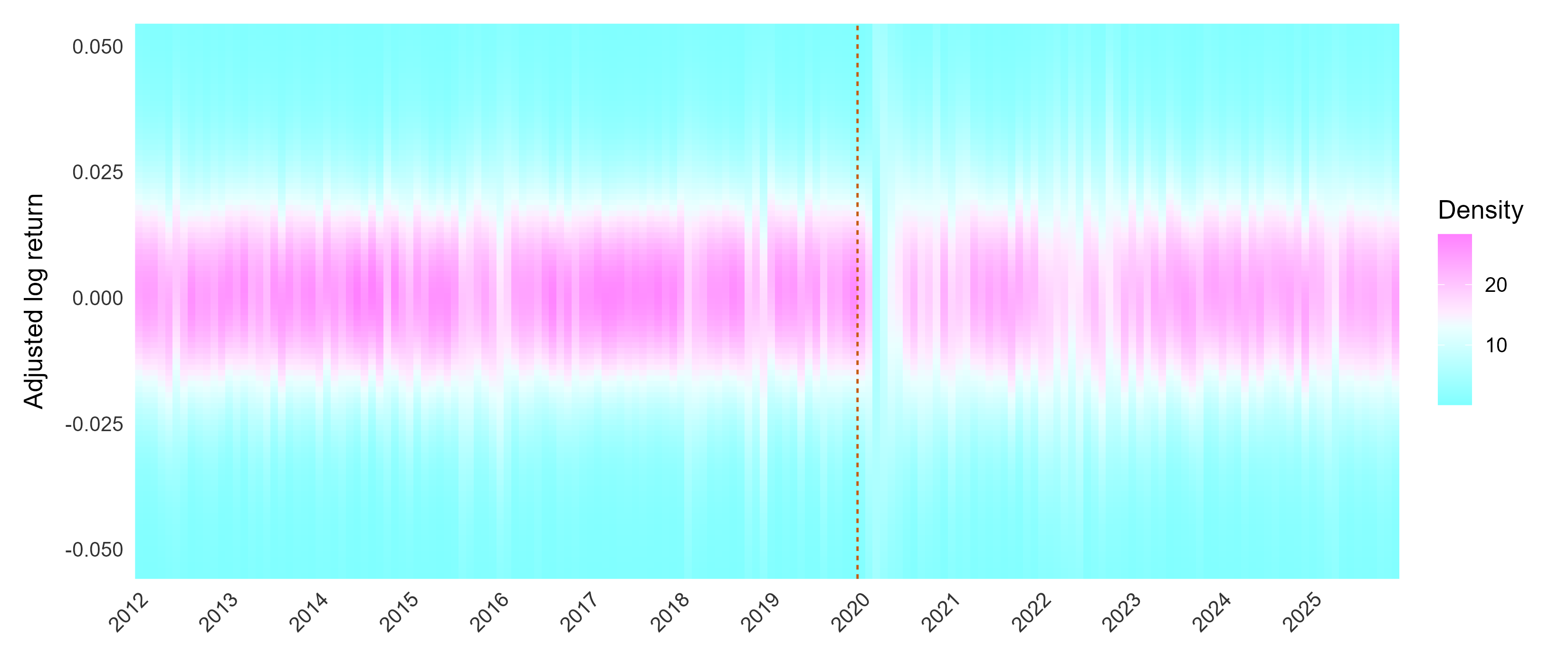}
		\caption{Estimated monthly distributions of pooled adjusted log returns from January 2012 to December 2025. The dashed vertical line marks January 2020.}
		\label{fig_stock_density}
	\end{figure}

	Following Section~5.1 of \citet{zhang2025change}, we estimate each monthly return distribution $X_t$ by a Gaussian kernel density estimator. Let $\widehat F_t$ denote the resulting estimated distribution function. We approximate the $2$-Wasserstein distance using 1001 deterministic interior quantiles:
	\begin{align*}
		\widehat W_2(X_i,X_j)
		={}&\left[
		\frac{1}{1001}\sum_{\ell=1}^{1001}
		\left\{\widehat F_i^{-1}(p_\ell)
		-\widehat F_j^{-1}(p_\ell)\right\}^2
		\right]^{1/2},
		\qquad
		p_\ell=\frac{\ell-1/2}{1001}.
	\end{align*}
	We then set $d(X_i,X_j)=\widehat W_2(X_i,X_j)$ and $\cK(Y_i,Y_j)=\exp\{-d(X_i,X_j)\}$. Figure~\ref{fig_stock_density} displays the estimated monthly return distributions. The dispersion of the density increases sharply in early 2020, as indicated by the vertical dashed line.
	
\begin{figure}[H]
		\centering
		\includegraphics[width=0.92\linewidth]{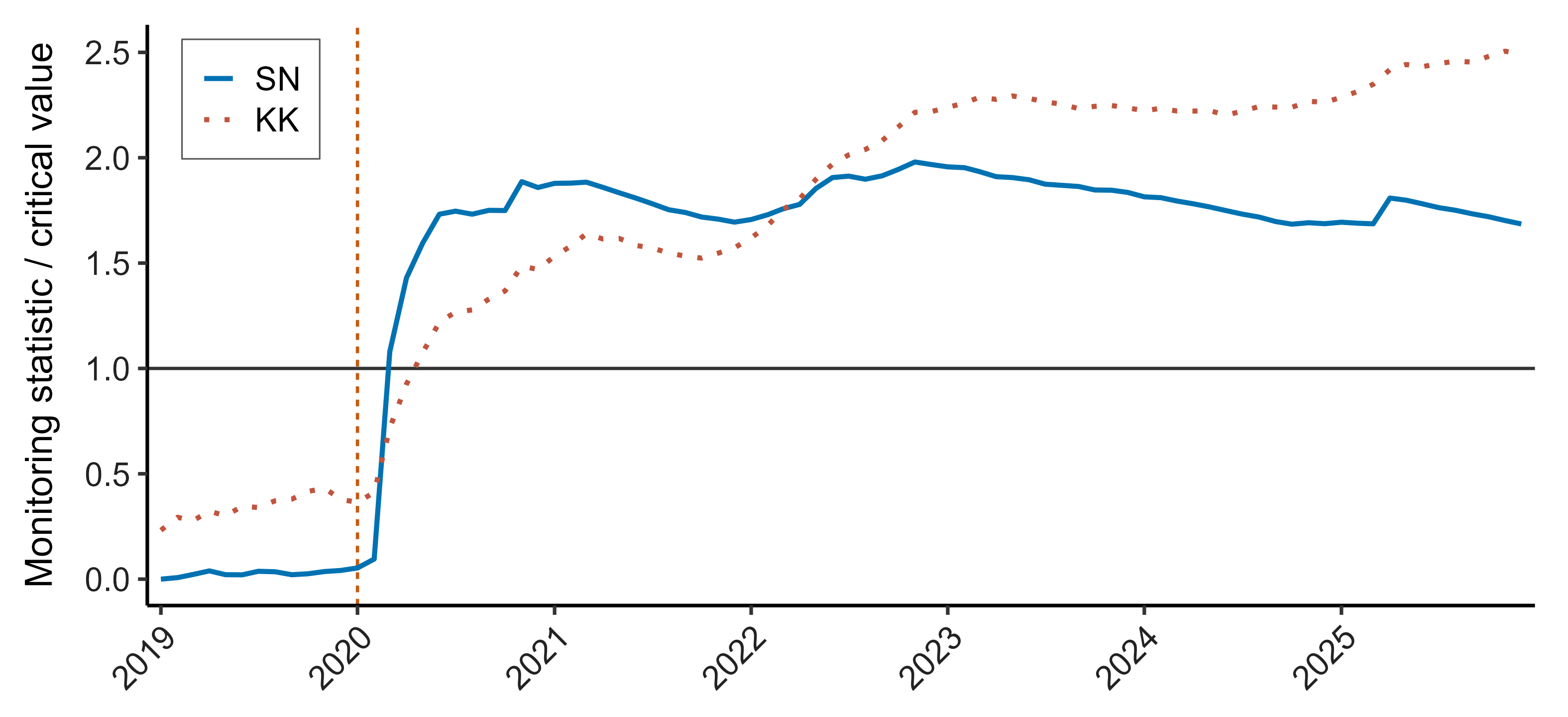}
		\caption{Monitoring paths for SN and KK with $n=84$, divided by their corresponding $95\%$ critical values. The horizontal line marks the critical threshold, and the vertical dashed line marks January 2020.}
		\label{fig_stock_monitor}
	\end{figure}

\begin{figure}[H]
		\centering
		\includegraphics[width=0.92\linewidth]{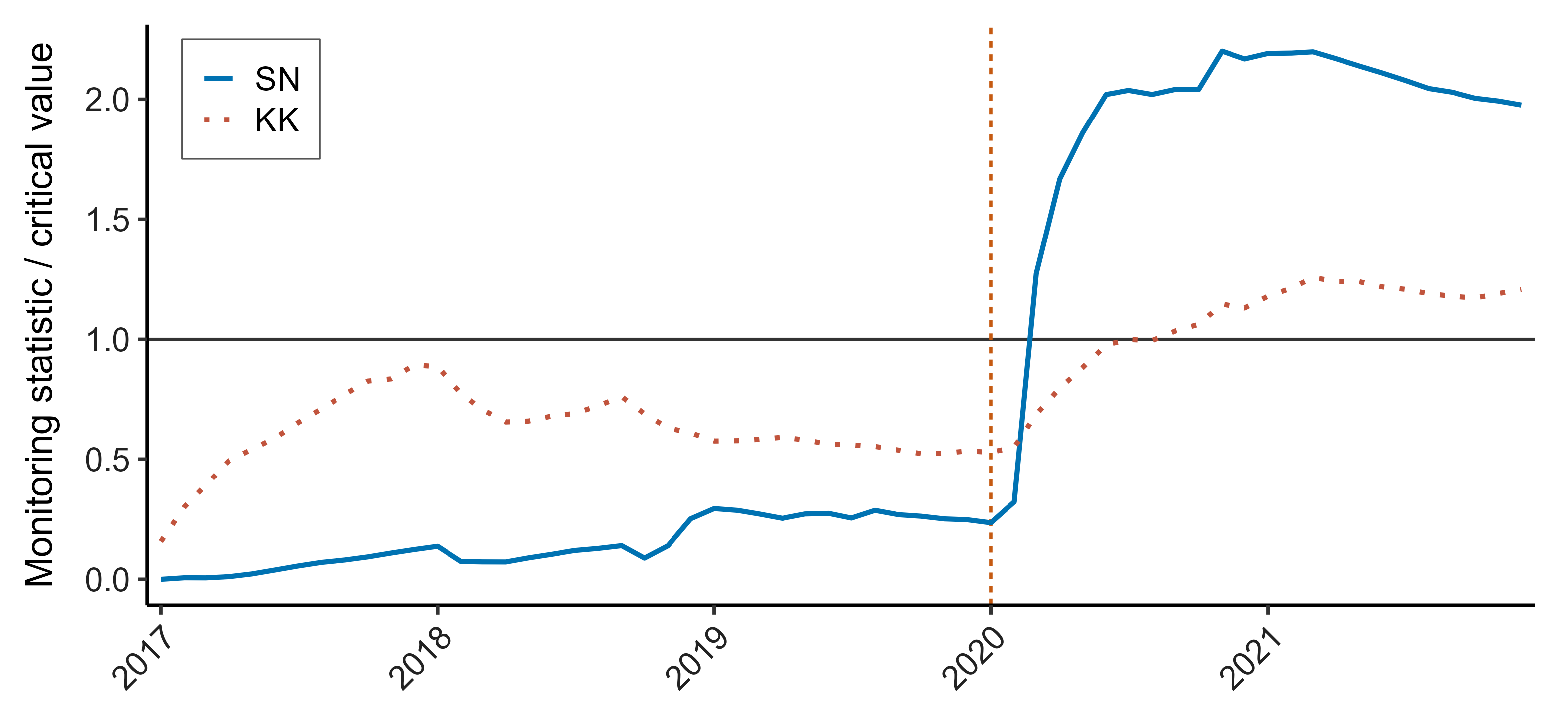}
		\caption{Monitoring paths for SN and KK with $n=60$, divided by their corresponding $95\%$ critical values. The horizontal line marks the critical threshold, and the vertical dashed line marks January 2020.}
		\label{fig_stock_monitor_n60}
	\end{figure}

	We apply SN, SN-PE, and the KK procedure of \citet{kutta2025monitoring}, using the specifications described in Section~\ref{sec3simu}.
	Figures~\ref{fig_stock_monitor} and~\ref{fig_stock_monitor_n60} display the SN and KK monitoring paths divided by their corresponding $95\%$ critical values, so the first crossing of one is the alarm time. Both historical-sample specifications lead to the same substantive conclusion for the proposed procedure: SN signals in March 2020, during the pronounced widening of the monthly return distributions visible in Figure~\ref{fig_stock_density}. Its alarm is unchanged when the historical period is shortened from seven to five years, indicating that the finding is not driven by the particular reference period used here. The power enhancement term in SN-PE vanishes in both analyses, so SN-PE reproduces the SN monitoring path and is omitted from the figures. KK also detects the early-2020 episode, but it signals later, in May 2020 with the seven-year historical sample and in September 2020 with the five-year historical sample. Thus, in this application, SN responds sooner and its detection time is more stable across the two historical-sample specifications, whereas the delay of KK is more sensitive to the choice of historical period.
	
	\section{Conclusion}\label{sec5conclu}
	
	In this paper, we propose closed- and open-end procedures for monitoring changes in the marginal distribution of a broad class of weakly dependent object-valued time series. By combining two distance-based Hilbert-space embeddings, a monitoring-time-dependent projection, and historical self-normalization, the detector can be computed entirely from pairwise distances and does not require long-run variance estimation. For a fixed horizon, it has a pivotal Brownian-functional null limit and admits an optional split-location weight. For an unbounded horizon, a growing monitoring-time weight controls the remote tail; a mixing maximal inequality and a converging-together argument yield a pivotal open-end limit, and Brownian time inversion gives a fixed-interval representation for critical-value simulation. Both procedures are omnibus-consistent against fixed marginal distributional changes and have local detection boundaries $n^{-1/2}$ when the first-order projection signal is nonzero and $n^{-1/4}$ when it vanishes and the quadratic signal becomes leading. Power enhancement has been shown to improve local power in the latter case.
	The detector also admits an exact iterative implementation. The simulations for distribution- and graph-valued time series in Section~\ref{sec3simu} and the supplementary material, respectively, support
	the closed-end calibration and show increasing power as the change
	magnitude grows, including under temporal dependence.  In the application to monthly stock-return distributions, SN produces a first alarm in March 2020 under both historical-sample specifications.
	
	The open-end alternative theory treats changes at $k^\ast=\lfloor nr_0\rfloor$ for fixed $r_0>1$. A distinct problem is to allow the scaled change location itself to diverge, $r_0=r_{0,n}\to\infty$; in that regime, detectability depends jointly on the change magnitude and the growth of the weight $g(r_{0,n})$. Another direction is to optimize the closed-end split-location weight and the open-end monitoring weight for detection delay subject to false-alarm control.

	{\fontsize{11.15}{11.15}\selectfont
\setlength{\parskip}{0pt}
\setlength{\bibsep}{0pt plus 0.24ex}

\putbib[cp_monitoring]
}
\end{bibunit}

\end{document}